\documentclass[12pt]{amsart}
\usepackage{color}
\usepackage{amsfonts,amsmath,amssymb}
\usepackage{latexsym,graphicx,xcolor}
\usepackage{txfonts}
 \usepackage{chicago}
\usepackage{fourier-orns}
\usepackage{multirow,rotating}
\usepackage{dsfont}
\usepackage{textcomp}
\newtheorem{defi}{Definition}[section]

\newtheorem{prop}{Proposition}[section]
\newtheorem{theo}{Theorem}[section]
\newtheorem{lem}{ Lemma}[section]

\newtheorem{remark}{Remark}[section]

\newcommand{\eref}[1]{(\ref{#1})}
\newcommand\ind{{{{1}}\hspace{-1,1mm}{\mathrm I}}}

\def\1{{\bf 1}}

\newcounter{hypc}
\newcommand{\hyp}[1]{\stepcounter{hypc}\tag{$\mathbf{A_{\thehypc}}$}
\label{#1}}

\def\build#1_#2^#3{\mathrel{\mathop{\kern 0pt#1}\limits_{#2}^{#3}}}

\def\cvL1{\build{\ \longrightarrow\ }_{\nti}^{{\mathbb{L}^1}}}

\def\nti{n \to \infty}

\title{Model selection in logistic regression}

 \author{  Marius Kwemou$^{(1)}$, Marie-Luce Taupin$^{(1)}$$^{(2)}$, Anne-Sophie Tocquet$^{(1)}$ }

 \address{(1) Laboratoire de Math\'ematiques et de Mod\'elisation d'Evry\\
 Universit\'e d'\'Evry Val d'Essonne\\
 UMR CNRS 8071- USC INRA \\
 23 Boulevard de France\\
 91037 \'Evry}

 \address{(2)
 INRA, UR 341 MIA-Jouy,\\
  Domaine de Vilvert,\\
  F78352 Jouy-en-Josas, France}

\begin{document}

\begin{abstract}
This paper is devoted to model selection in logistic regression.  We   extend the model selection principle introduced by  Birg\'e and Massart \citeyear{birge2001gaussian} to logistic regression model. This selection is done by using penalized maximum  likelihood criteria.
We propose in this context a completely data-driven criteria based on the  slope heuristics. We prove  non asymptotic oracle inequalities for  selected estimators.
Theoretical results are illustrated through  simulation studies.
\end{abstract}

\maketitle



\medskip


\noindent {\bf AMS 2000 MSC}: Primary 62J02, 62F12, Secondary 62G05, 62G20.
\section{Introduction}
\setcounter{equation}{0}
\setcounter{lem}{0}
\setcounter{theo}{0}

 Consider the following generalization of the logistic regression model : let  $(Y_1,x_1),\cdots,(Y_n,x_n)$, be a sample of size $n$  such that
$(Y_i, x_i)\in \{0,1\}\times\mathcal{X}$ and 
\begin{eqnarray*}
\mathbb{E}_{f_0}(Y_i)=\pi_{f_0}(x_i)=\frac{\exp{f_{0}(x_i)}}{1+\exp{f_{0}(x_i)}},
\end{eqnarray*}
where  $f_{0}$ is an unknown function to be estimated and the design points $x_1,...,x_n$ are deterministic. This model can be viewed as a nonparametric version of the "classical" logistic model which  relies on the assumption that $x_i\in \mathbb{R}^d$, and that there exists $\beta_{0}\in \mathbb{R}^d$ such that $f_0(x_i)=\beta_0^\top x_i.$

Logistic regression is a widely used model  for predicting the outcome of binary dependent variable.  For example  logistic model can be used in medical study to predict the probability that a patient has a given disease (e.g. cancer), using observed characteristics (explanatory  variables) of the patient such as weight, age, patient's gender \textit{etc}. However  in the presence of numerous  explanatory variables  with potential influence,  one would like to use only a few number of variables, for the sake of interpretability or to avoid overfitting.  But it is not always obvious to choose the adequate variables. This is the well-known problem of variables selection or model selection. 

In this paper, the unknown function  $f_0$ is not specified and not necessarily linear. Our aim is to estimate $f_0$  by a  linear combination  of  given functions, often called dictionary. The dictionary can be a basis of functions, for  instance spline or polynomial basis.  

A nonparametric version of the classical logistic model  has already been considered by Hastie \citeyear{non_par_logist}, where a nonparametric estimator of $f_0$ is proposed using local maximum likelihood.
The problem of nonparametric estimation in additive regression model is well known and deeply studied.  But in logistic regression model it is less studied.  One can cite for instance Lu \citeyear{Lu2006}, Vexler \citeyear{Vexler2006}, Fan \textit{ et al.} \citeyear{Fanetal1998}, Farmen \citeyear{Farmen1996},  Raghavan \citeyear{Raghavan1993}, and Cox \citeyear{Cox1990}.

Recently few papers deal with model selection or nonparametric estimation in logistic regression using  $\ell_1$ penalized contrast Bunea \citeyear{bunea2008honest}, Bach \citeyear{bach10}, van de Geer \citeyear{van2008}, Kwemou \citeyear{kwemou}. Among them, some establish non asymptotic oracle inequalities that hold even in high dimensional setting.
When the dimension of $\mathcal{X}$ is high, that is greater than dozen, such $\ell_1$ penalized contrast estimators are known to provide reasonably good results.
When the dimension of $\mathcal{X}$ is small,  it is often better to choose different penalty functions. One classical penalty function is what   we call $\ell_0$ penalization. 
Such penalty functions, built as increasing function of the dimension of $\mathcal{X}$, usually refers to model selection. The last decades have witnessed a growing interest in the model selection problem since the seminal works of Akaike~\citeyear{akaike1973}, Schwarz \citeyear{schwarz1978}.  In additive regression one can cite among the others  Baraud~\citeyear{baraud2000model}, Birg\'e and Massart \citeyear{birge2001gaussian}, Yang \citeyear{yang1999}, in density estimation Birg\'e~\citeyear{birge2014model}, Castellan \citeyear{castellan2003density} and in  segmentation problem Lebarbier \citeyear{lebarbier}, Durot \textit{et al.} \citeyear{DurotLebarbierTocquet}, and Braun \textit{ et al.} \citeyear{braun}. All the previously cited papers use $\ell_{0}$ penalized contrast to perform model selection. But model selection procedures based on penalized maximum likelihood estimators in logistic regression are less studied in the literature.

 In this paper we focus on  model selection using $\ell_0$ penalized contrast for logistic regression model and in this context we state non asymptotic oracle inequalities. More precisely, given some collection functions, we consider estimators of $f_0$ built as linear combination of the functions.
The point that the true function is not supposed to be  linear combination of those  functions, but we expect that the spaces of linear combination of those functions would provide suitable approximation spaces. Thus, to this collection of functions, we associate a collection of estimators of $f_0$.
Our aim is to propose a data driven procedure, based on penalized criterion, which will be able to choose
the "best" estimator among the collection of estimators, using $\ell_{0}$ penalty functions. 
 
The collection of estimators is built using minimisation of the opposite of logarithm likelihood. The properties of estimators are described in term
of Kullback-Leibler divergence  and the empirical $L_2$ norm.  Our results can be splitted into two parts.

 First, in a general model selection framework, with general collection of functions we provide a completely data driven procedure
that automatically selects the best model among the collection. We state  non asymptotic oracle  inequalities  for Kullback-Leibler divergence  and the empirical $L_2$ norm  between the selected estimator and the true function $f_0$. The estimation procedure relies on the building of a suitable penalty function, suitable in the sense that it performs best risks and suitable in the sense
that it does not depend on the unknown smoothness parameters of the true function $f_0$. But, 
the  penalty function depends on a bound related to target function $f_0$. This can be seen as the price to pay  for the generality. It comes from
needed links between   Kullback-Leibler divergence and empirical $L_2$ norm. 

Second, we consider the specific case of collection of piecewise functions  which provide estimator of type regressogram. In this case,  we exhibit a completely data driven penalty, free from $f_0$. The model selection procedure based on this penalty provides an  adaptive estimator and state a non asymptotic oracle inequality for Hellinger distance  and the empirical $L_2$ norm  between the selected estimator and the true function $f_0$.
In the case of piecewise constant functions basis, the connection between Kullback-Leibler divergence  and the empirical $L_2$ norm are obtained without bound on the true function $f_0$.
This last result is of great interest for example in segmentation study, where the target function is piecewise constant or can be well approximated by piecewise constant functions.

Those theoretical results are illustrated through simulation studies. In particular we show that our  model selection procedure (with the suitable penalty) have good non asymptotic properties
as compared to usual known criteria such as AIC and BIC. A great attention has been made on the practical calibration of the penalty function. This practical calibration is mainly based on the ideas
of what is usually referred as slope heuristic as proposed in Birg\'e and Massart \citeyear{birge2007} and developed in Arlot and Massart \citeyear{arlot2009}.

The  paper is organized as follow. In Section~\ref{S1} we set our framework and describe our estimation procedure. In Section~\ref{S2} we  define  the model selection procedure  and state the oracle inequalities in the general framework. Section~\ref{S3} is devoted to regressogram  selection, in this section, we establish a bound of the Hellinger risk between the selected model and the target function. The simulation study is reported in Section~\ref{S4}. The proofs of the results are postponed to Section~\ref{S5} and ~\ref{S6}.  
\section{Model and framework}\label{S1}

\setcounter{equation}{0}
\setcounter{lem}{0}
\setcounter{theo}{0}
 Let  $(Y_1,x_1),\cdots,(Y_n,x_n)$, be a sample of size $n$  such that
$(Y_i, x_i)\in \{0,1\}\times\mathcal{X}$. Throughout the paper, we consider a fixed design setting \textit{i.e.}  $x_{1},\dots,x_{n}$ are considered as deterministic. In this setting, consider the extension of the "classical" logistic regression model \eref{model} where we aim at estimating the unknown function $f_{0}$
in
\begin{eqnarray}
\label{model}
\mathbb{E}_{f_0}(Y_i)=\pi_{f_0}(x_i)=\frac{\exp{f_{0}(x_i)}}{1+\exp{f_{0}(x_i)}}.
\end{eqnarray}
We propose to estimate the unknown function $f_0$ by model selection. This model selection is performed using penalized maximum likelihood estimators. In the following we denote by $\mathbb{P}_{f_0}(x_1)$  the distribution of $Y_1$ and by $\mathbb{P}^{(n)}_{f_0}(x_1,\cdots,x_n)$ the distribution of   $(Y_1,\dots,Y_n)$ under Model (\ref{model}). Since the variables $Y_i$'s are independent random variables,
$$\mathbb{P}^{(n)}_{f_0}(x_1,\cdots,x_n)=\Pi_{i=1}^n \mathbb{P}_{f_0}(x_i)  =\prod_{i=1}^n \pi_{f_0}(x_i)^{Y_i}(1-\pi_{f_0}(x_i))^{1-Y_i}  .$$
It follows that  for a function $f$ mapping $\mathcal{X}$ into $\mathbb{R}$, the likelihood is defined as:
\begin{eqnarray*}
L_n(f) = \mathbb{P}^{(n)}_{f}(x_1,\cdots,x_n)  =\prod_{i=1}^n \pi_f(x_i)^{Y_i}(1-\pi_f(x_i))^{1-Y_i},
\end{eqnarray*}
where
\begin{eqnarray}
\label{pif}
\pi_{f}(x_i)=\frac{\exp{(f(x_i))}}{1+\exp({f}(x_i))}.
\end{eqnarray}
We choose the opposite of the log-likelihood as the  estimation criterion that is
\begin{eqnarray}\label{gamma_n}
\gamma_n(f)=-\frac{1}{n}\log(L_n(f))=\frac{1}{n}\sum_{i=1}^{n}\Big\{\log(1+e^{f(x_i)})-Y_if(x_i)\Big\}.
\end{eqnarray}
Associated to this estimation criterion we consider the Kullback-Leibler information divergence
 $\mathcal{K}(\mathbb{P}_{f_0}^{(n)},\mathbb{P}_{f}^{(n)})$ defined
as
\begin{eqnarray*}
 \mathcal{K}(\mathbb{P}_{f_0}^{(n)},\mathbb{P}_{f} ^{(n)})=\frac{1}{n}\int \log\left( \frac{\mathbb{P}_{f_0}^{(n)}}{\mathbb{P}_{f}^{(n)}}\right) d\mathbb{P}_{f_0}^{(n)}.
\end{eqnarray*}
The loss function  is the excess risk, defined as
 \begin{eqnarray}
\label{gamma}\mathcal{E}(f):=\gamma(f)-\gamma(f_0)  \mbox{ where, for any }f, \quad
\gamma(f)=\mathbb{E}_{f_0}[\gamma_n(f)].~~~~~
\end{eqnarray}
Easy calculations show that the excess risk is linked to the Kullback-Leibler information divergence through the relation
 $$\mathcal{E}(f)=\gamma(f)-\gamma(f_0)=\mathcal{K}(\mathbb{P}_{f_0}^{(n)},\mathbb{P}_{f}^{(n)}).$$
It follows that, $f_0$ minimizes the excess risk, that is 
$$f_0= \arg\min_{f} \gamma (f).$$ 
As usual, one can not  estimate $f_0$ by the minimizer of $\gamma_n(f)$ over any functions space, since it is infinite.
The usual way is to minimize $\gamma_n(f)$ over a finite dimensional collections of models, associated to a finite dictionary of functions  $\phi_j : \mathcal{X} \rightarrow \mathbb{R}$
$$\mathcal{D}=\{\phi_1,\dots,\phi_M\}.$$
For the sake of simplicity  we will suppose that  $\mathcal{D}$ is a orthonormal basis of functions. Indeed,  if  $\mathcal{D}$ is not an orthonormal basis of functions, we can always find an orthonormal basis
of functions  $\mathcal{D^{\prime}}=\{\psi_1,\dots,\psi_{M^{\prime}}\}$ such that 
$$\langle \phi_1,\dots,\phi_M \rangle=\langle \psi_1,\dots,\psi_{M^{\prime}} \rangle.$$

Let  $\mathcal{M}$ the set of all subsets  $m \subset \{1,\dots,M\}$.  For every $m\in \mathcal{M}$,  we call $\mathcal{S}_m$ the model  
\begin{equation}\label{model_col_gen}
\mathcal{S}_m:=\Big\{f_{\beta}=\sum_{j\in m}\beta_j\phi_j\Big \}
\end{equation}
 and  $D_m$ the dimension of the span of $\{\phi_j, j \in m\}$.
Given the countable collection of models $\{\mathcal{S}_{m}\}_{m\in \mathcal{M}}$,  we define $\{\hat f_m\}_{m\in \mathcal{M}}$ the corresponding estimators, \textit{i.e.} the estimators obtaining by minimizing  $\gamma_n$ over each model $\mathcal{S}_m$. 
For each $m\in \mathcal{M}$, $\hat f_m$ is defined by
\begin{equation}\label{fchapD1}
 \hat f_m = \arg\min_{t\in \mathcal{S}_m} \gamma_n(t).
\end{equation}

Our aim is choose the "best" estimator among this collection of estimators, in the sense that 
it minimizes the risk. In many cases, it is not easy  to choose  the "best" model.  Indeed,  a model with small dimension tends to be efficient from estimation point of view whereas it could be far from the "true" model. On the other side, a more complex model easily fits data but the estimates have poor predictive performance (overfitting). We thus expect that this best estimator mimics  what is usually called the oracle defined as
\begin{equation}\label{m*}
m^{*}=\arg\min_{m\in \mathcal{M}} \mathcal{K}(\mathbb{P}_{f_0}^{(n)},\mathbb{P}_{\hat f_m}^{(n)}).
\end{equation}
Unfortunately,  both, minimizing the risk and minimazing the kulback-leibler divergence, require the knowledge of  the true (unknown) function $f_0$ to be estimated.

Our goal is to  develop a data driven strategy based on data, that automatically selects the best estimator among the collection, this best estimator having a risk as close as possible to the oracle risk, that is the risk of $\hat f_{m^{*}}$. In this context, our strategy follows the lines of  model selection as developed by Birg\'e and Massart \citeyear{birge2001gaussian}. We also refer to the book
Massart \citeyear{massart2007} for further details on model selection.

 We use penalized maximum likelihood estimator for choosing some data-dependent $\hat m$ nearly as good as the ideal choice $m^{*}$.
More precisely, the idea is to select $\hat{m}$ as a minimizer of the penalized criterion 
\begin{equation}\label{crit0}
\hat{m}=\arg \min_{m\in  \mathcal{M}}\left\{  \gamma_n(\hat f_m)+\mbox{pen}(m)\right\},
\end{equation}
where $\mbox{pen} : \mathcal{M} \longrightarrow \mathbb{R}^{+}$ is a  data driven penalty function.
The estimation properties of $\hat f_m$  are evaluated by non asymptotic bounds of a risk associated to a suitable chosen loss function.
 The great challenge is choosing  the penalty  function such that the selected model $\hat m$ is 
nearly as good as the oracle $m{^*}$. This penalty term is classically based on the idea that
$$m^{*}=\arg\min_{m\in \mathcal{M}} \mathbb{E}_{f_0}\mathcal{K}(\mathbb{P}_{f_0}^{(n)},\mathbb{P}_{\hat f_m}^{(n)})=\arg\min_{m \in \mathcal{M}} \left[ \mathbb{E}_{f_0}\mathcal{K}(\mathbb{P}_{f_0}^{(n)},\mathbb{P}_{f_m}^{(n)})  + \mathbb{E}_{f_0}\mathcal{K}(\mathbb{P}_{f_m}^{(n)},\mathbb{P}_{\hat f_m}^{(n)})  \right]$$
where $f_m$ is defined as
\begin{equation*}
f_m = \arg\min_{t \in S_m} \gamma(t).
\end{equation*}
Our goal is to build a penalty function such that the selected model $\hat m$ fulfills an oracle inequality:
$$\mathcal{K}(\mathbb{P}_{f_0}^{(n)},\mathbb{P}_{\hat f_{\hat m}}^{(n)})\leq C_n\inf_{m\in \mathcal{M}}\mathcal{K}(\mathbb{P}_{f_0}^{(n)},\mathbb{P}_{\hat f_m}^{(n)}) + R_n.$$
This inequality is expected to hold either in expectation or with high probability, where $C_n$ is as close to 1 as possible and $R_n$ is a remainder term negligible  compared to  $\mathcal{K}(\mathbb{P}_{f_0}^{(n)},\mathbb{P}_{\hat f_{m{*}}}^{(n)})$. 

In the following we consider two separated case. 
First we consider  general collection of models under boundedness assumption. Second we consider the specific case of regressogram collection.
\section{Oracle inequality for general models collection under boundedness assumption}
\label{S2}
Consider model (\ref{model}) and 
$(\mathcal{S}_m)_{m\in\mathcal{M}}$ a collection of models defined by (\ref{model_col_gen}). 
Let $C_0 >0$ and $\mathbb{L}_\infty(C_0)=\Big\{f : \mathcal{X}\to \mathbb{R}, ~~ \max_{1\leqslant i \leqslant n}|f(x_{i})|\leqslant C_0\Big\}$. For $m\in\mathcal{M}$,  $\gamma_n$ given in (\ref{gamma_n}), and $\gamma$ is given by (\ref{gamma}), we define  
\begin{equation}\label{fchapD1c}
 \hat f_m = \arg\min_{t\in \mathcal{S}_m\cap \mathbb{L}_\infty(C_0)} \gamma_n(t) 
\mbox{ and }
f_m = \arg\min_{t \in S_m\cap \mathbb{L}_\infty(C_0)} \gamma(t).
\end{equation}

The first step consists in studying the estimation properties of $\hat f_m$ for each $m$, as it is stated in the following proposition. 
\begin{prop}\label{borne1}
Let $C_0 >0$ and  $\mathcal{U}_0=e^{C_0}/(1+e^{C_0})^2$. For $m\in\mathcal{M}$, let $ \hat f_m$ and $f_m$ as in (\ref{fchapD1c}).
 We have
\begin{eqnarray*}
\mathbb{E}_{f_0} [\mathcal{K}(\mathbb{P}_{f_0}^{(n)},\mathbb{P}_{\hat{f}_m}^{(n)})]\leqslant \mathcal{K}(\mathbb{P}_{f_0}^{(n)},\mathbb{P}_{f_m}^{(n)})+\frac{D_m}{2 n \mathcal{U}_0^{2}}
\end{eqnarray*}
\end{prop}
\noindent
This proposition says that  the "best" estimator amoung the collection $\{\hat f_m\}_{m\in\mathcal{M}}$, in the sense of the Kullback-Leibler risk,   is the one  
which makes a balance between the bias and the complexity of the model.
 In the ideal situation where  $f_0$ belongs to $\mathcal{S}_m$, we have that
\begin{eqnarray*}
\mathbb{E}_{f_0} [\mathcal{K}(\mathbb{P}_{f_0}^{(n)},\mathbb{P}_{\hat{f}_m}^{(n)})]\leqslant \frac{1}{\mathcal{U}_0^{2}}\frac{D_m}{2n}. 
\end{eqnarray*}

To derive the model selection procedure we need the following assumption : 
\begin{align}
& \hyp{A1}  \mbox{There exists a constant}~  0 < c_{1} < \infty~ \mbox{such that } ~\max_{1\leqslant i \leqslant n}|f_{0}(x_{i})|\leqslant c_{1}.~~~~~~~~~~~~~~~~~~~~~~~~~~~~~~~~~~~
\end{align}
In the following theorem we propose a choice for the penalty function and we state non asymptotic risk bounds.  \begin{theo}\label{theo1} 
Given $C_0 >0$, for $m\in\mathcal{M}$, let $ \hat f_m$ and $f_m$ be defined as (\ref{fchapD1c}).
Let us denote  $\parallel f\parallel_n^2=n^{(-1)}\sum_{i=1}^n f^2(x_i)$.
 Let  $\{L_m\}_{m\in \mathcal{M}}$ some positive numbers  satisfying 
 $$\Sigma=\sum_{m\in \mathcal{M}}\exp(-L_mD_m)< \infty.$$
 We define $\mbox{pen}: \mathcal{M} \rightarrow \mathbb{R}_{+}$ , such that,  for $m\in \mathcal{M}$, 
$$\mbox{pen}(m)\geqslant \lambda\frac{D_m}{n}\left(\frac{1}{2}+\sqrt{5L_m}\right)^{2}, $$
where $\lambda$ is a positive constant depending on $c_1$. Under Assumption~(\ref{A1}) we have
$$\mathbb{E}_{f_0}[\mathcal{K}(\mathbb{P}_{f_0}^{(n)},\mathbb{P}_{\hat f_{\hat m}}^{(n)})]\leqslant C\inf_{m\in \mathcal{M}}\left\{\mathcal{K}(\mathbb{P}_{f_0}^{(n)},\mathbb{P}_{f_m}^{(n)})+\mbox{pen}(m)\right\}+C_1\frac{\Sigma}{n}$$
and
$$\mathbb{E}_{f_0}\parallel \hat f_{\hat{m}}-f_{0}\parallel_{n}^{2}\leqslant C^{\prime}\inf_{m\in\mathcal{M}}\left\{ \parallel f_{0}-f_m\parallel_{n}^{2}+\mbox{pen}(m)\right\}+C_{1}^{\prime}\frac{\Sigma}{n}.$$
where $C,C^{\prime}, C_{1},C_{1}^{\prime}$ are constants depending on $c_{1}$ and $C_{0}$. 
\end{theo}
\noindent
This theorem provides oracle inequalities for $L_{2}-$norm and for K-L divergence between the selected model and the true function.
Provided that penalty has been properly chosen, one can bound the   $L_{2}-$norm and the K-L divergence between the selected model and the true function. The inequalities in Theorem~\ref{theo1} are non-asymptotic inequalities in the sense that the result is obtain for a fixed $n$. This theorem is very general and does not make specific assumption on the dictionary. However,  the penalty function depends on some unknown constant  $\lambda$ which depends on the bound of the true function $f_0$ through Condition \eref{cond_lambda}. In practice this constant can be calibrated using "slope heuristics" proposed in Birg\'e and Massart \citeyear{birge2007}. In the following we will show how to obtain similar result with a penalty function not connected to the bound of the true unknown function $f_0$ in the regressogram case. 
\section{Regressogram  functions}\label{S3}
\subsection{Collection of models }\label{models}
In this section we suppose (without loss of generality) that $f_0 : [0,1]\to \mathbb{R}$.  
For the sake of simplicity, we use the notation  $f_0(x_i)=f_0(i)$ for every $i=1,\dots,n$. Hence  $f_0$ is defined from $\{1,\dots,n\}$ to $\mathbb{R}$.
Let $\mathcal{M}$ be a collection of partitions of intervals of $\mathcal{X}=\{1,\dots,n\}$. For any $m\in\mathcal{M}$ and $J\in m$, let $\ind_J$ denote the indicator function of $J$ and $S_m$ be the linear span of $\{\ind_J,J\in m\}$. 
When  all intervals have the same length, the partition is said regular, and is is irregular otherwise.

\subsection{Collection of estimators: regressogram}
For  a  fixed $m$, the minimizer  $\hat f_m$  of the
empirical contrast function $\gamma_n$, over $S_m$, is called the \textit{regressogram}.  That is, 
 $f_0$ is estimated by $\hat{f}_m$ given by
\begin{eqnarray}
\label{fchapD}
\hat{f}_m
=\arg\min_{f\in S_m}\gamma_n(f).
\end{eqnarray}
where $\gamma_n$ is given by (\ref{gamma_n}). Associated to $S_m$ we have
\begin{eqnarray}
\label{fD} f_m=\arg\min_{f\in S_m}\gamma(f)-\gamma(f_0)=\arg\min_{f\in S_m}\mathcal{K}(\mathbb{P}_{f_0}^{(n)},\mathbb{P}_{f}^{(n)}).
\end{eqnarray}
In the specific case where $S_m$  is the set of piecewise constant functions on some partition $m$,  $\hat{f}_m$ and ${f}_m$ are given by the following lemma. 
\begin{lem} \label{Projhisto} 
For $m\in \mathcal{M}$ , let $f_m$ and $\hat f_m$ be defined by \eref{fD} and \eref{fchapD} respectively . 
Then, $f_m=\sum_{J\in m}\overline{f}_m^{(J)}\ind_{J}$ and  $\hat f_m=\sum_{J\in m}\hat{f}_m^{(J)}\ind_{J}$ with
$$\overline{f}_m^{(J)}=\log \left(  
\frac{ \sum_{i\in J}\pi_{f_0}(x_i)  }{\vert  J \vert  (1-\sum_{i\in J}\pi_{f_0}(x_i) / \vert  J \vert ) }\right)
\mbox{ and } \hat{f}_m^{(J)}=\log \left(  
\frac{ \sum_{i\in J}Y_i  }{\vert  J \vert  (1-\sum_{i\in J}Y_i/ \vert  J \vert ) }\right).$$
Moreover,  $\pi_{f_m}=\sum_{J\in m}\pi_{f_m}^{(J)}\ind_{J}$ and  $\pi_{\hat f_m}=\sum_{J\in m}\pi_{\hat f_m}^{(J)}\ind_{J}$ with
$$\pi_{f_m}^{(J)}=\frac{1}{\vert J\vert}\sum_{i\in J}\pi_{f_0}(x_i), \mbox{ and } \pi_{\hat f_m}^{(J)}=\frac{1}{\vert J\vert}\sum_{i\in J}Y_i.$$
\end{lem}
Consequently, $\pi_{f_m}=\arg\min_{\pi \in S_m}\parallel \pi-\pi_{f_0}\parallel^2_n$
is the usual projection of $\pi_{f_0}$ on to $S_m$.


\subsection{First bounds on $\hat{f}_m$}
\setcounter{equation}{0}
\setcounter{lem}{0}
\setcounter{theo}{0}
Consider the following assumptions:
\begin{align}
&\hyp{H0} \mbox{ There exists a constant }\rho>0 \mbox{ such that }
\min_{i=1,\cdots,n}\pi_{f_0}(x_i)\geq \rho~~ \mbox{and}~~~ \min_{i=1,\cdots,n} [1-\pi_{f_0}(x_i)]\geq \rho .
\end{align}

\begin{prop}\label{borne2}
Consider Model \eref{model} and let $\hat{f}_m$ be defined by \eref{fchapD} with $m$ such that for all $J\in m$, $\vert J\vert \geqslant \Gamma [\log(n)]^2$ for a positive constant $\Gamma$. Under Assumption~(\ref{H0}),  for all $\delta>0$ and $a>1$, we have 
\begin{eqnarray*}
\mathbb{E}_{f_0}[\mathcal{K}(\mathbb{P}_{f_0}^{(n)},\mathbb{P}_{\hat f_{m}}^{(n)})]
&\leqslant& \mathcal{K}(\mathbb{P}_{f_0}^{(n)},\mathbb{P}_{f_{m}}^{(n)}))+\frac{(1+\delta)D_m}{(1-\delta)^2n} +\frac{\kappa(\Gamma,\rho,\delta)}{n^a}.
\end{eqnarray*}
\end{prop}

\subsection{Adaptive estimation and oracle inequality}\label{modsel}
The following result  provides an adaptive estimation of $f_0$ and  a risk bound of the selected model. 

\begin{defi}\label{def_mf}
Let $\mathcal{M}$ be a collection of partitions of
$\mathcal{X}=\{1,\dots,n\}$ constructed on the partition $m_f$ \textit{i.e.} $m_f$ is a refinement of every $m \in \mathcal{M}.$
\end{defi}
In other words, a partition $m$ belongs to $\mathcal{M}$ if any element of $m$ is the union of some elements of $m_f$. Thus $S_{m_f}$ contains every model of the collection  $\{S_m\}_{m\in \mathcal{M}}$.
\setcounter{equation}{0}
\setcounter{lem}{0}
\setcounter{theo}{0}
\begin{theo}\label{theo2}
 Consider Model \eref{model} under Assumption~(\ref{H0}). Let $\{S_m, m\in \mathcal{M}\}$ be a collection of models defined in Section \ref{models} where $\mathcal{M}$ is a set of partitions constructed on the partition $m_{f}$  such that 
 \begin{equation}\label{def_Gam}
 \mbox{for all} ~~J\in m_{f}, \vert J\vert \geq \Gamma \log^{2}(n),
 \end{equation}
 where $\Gamma$ is a  positive constant.  Let $(L_m)_{m\in \mathcal{M}}$  be some family of positive weights 
 satisfying 
 \begin{equation}\label{sigma}
 \Sigma=\sum_{m\in \mathcal{M}}\exp(-L_m D_m) < +\infty.
 \end{equation}
 Let $\mbox{pen}: \mathcal{M} \rightarrow \mathbb{R}_{+}$ satisfying for $m\in \mathcal{M}$, and for $\mu > 1,$ 
$$\mbox{pen}(m)\geqslant \mu\frac{D_m }{n}\left(1+6L_m+8\sqrt{L_m}\right).$$
 Let $\tilde{f}=\hat{f}_{\hat m}$ where $$\hat m =\arg\min_{m \in \mathcal{M}}\left\{\gamma_n(\hat f_m)+\mbox{pen}(m)\right\},$$
then, for  $C_{\mu}= 2\mu^{1/3}/(\mu^{1/3}-1)$, we have
\begin{equation}\label{oracle}
\mathbb{E}_{f_0}[h^2(\mathbb{P}^{(n)}_{f_0},\mathbb{P}^{(n)}_{\tilde f})]\leqslant C_{\mu}\inf_{m\in \mathcal{M}}\left\{\mathcal{K}(\mathbb{P}_{f_0}^{(n)},\mathbb{P}_{f_m}^{(n)})+\mbox{pen}(m)\right\}+\frac{C(\rho,\mu, \Gamma,\Sigma)}{n}.
\end{equation}

\end{theo}

This theorem provides a non asymptotic bound for the Hellinger  risk between the selected model and the true one. On the opposite of  Theorem~\ref{theo1}, the penalty function does not depend on the bound of the true function. The selection procedure  based only on the data offers the advantage to free the estimator from any prior knowledge  about the smoothness of the function to estimate. The estimator is therefore  adaptive.    
As we bound Hellinger risk in ~\eref{oracle} by Kulback-Leibler risk, one should prefer to have the Hellinger risk on the right hand side  instead of the Kulback-Leibler risk. Such a bound is possible if we assume that $\log(\Vert \pi_{f_0}/\rho\Vert_{\infty})$ is bounded. Indeed if we assume that there exists $T$ such that $\log(\Vert \pi_{f_0}/\rho\Vert_{\infty})\leq T$, this implies that $\log(\Vert \pi_{f_0}/\pi_{f_m}\Vert_{\infty})\leq T$ uniformly for all partitions $m\in \mathcal{M}.$ Now using Inequality (7.6) p. 362 in Birg\'e and Massart \citeyear{birge1998} we have that $\mathcal{K}(\mathbb{P}_{f_0}^{(n)},\mathbb{P}_{f_m}^{(n)})\leq (4+2\log(M))h^2(\mathbb{P}_{f_0},\mathbb{P}_{f_m})$ which implies, 
\begin{equation*}
\mathbb{E}_{f_0}[h^2(\mathbb{P}_{f_0}^{(n)},\mathbb{P}_{\tilde f}^{(n)})]\leqslant C_{\mu}.C(T)\inf_{m\in \mathcal{M}}\left\{h^2(\mathbb{P}_{f_0}^{(n)},\mathbb{P}_{f_m}^{(n)})+\mbox{pen}(m)\right\}+\frac{C(\rho,\mu, \Gamma,\Sigma)}{n}.
\end{equation*}

\subsubsection*{Choice of the weights  $\{L_m, m\in \mathcal{M}\}$} According to Theorem~\ref{theo2}, the penalty function depends on the collection $\mathcal{M}$ through the choice of the weights $L_m$ satisfying \eref{sigma}, \textit{i.e.}
\begin{equation}\label{Sigma1}
\Sigma=\sum_{m\in \mathcal-{M}}\exp(-L_m D_m) =\sum_{D\geq 1}e^{-L_DD} Card\{m\in \mathcal{M}, \vert m \vert= D\}<  \infty.
\end{equation}

Hence the number of models having the same dimension $D$ plays an important role
in the risk bound.

If there is only one model of dimension $D$, 
a simple way  of choosing $L_D$ is to take them constant, \textit{i.e.} $L_D=L$ for all $m\in \mathcal{M}$, and thus we have from~\eref{Sigma1}
$$\Sigma=\sum_{D\geq 1} e^{-LD}< \infty.$$ 

This is the case when  $\mathcal{M}$ is  a family of regular partitions.
Consequently, the choice \textit{i.e.} $L_D=L$ for all $m\in \mathcal{M}$ leads to a penalty proportional to the dimension $D_m$, and for every $D_m\geq1$,
\begin{equation}\label{penlin}
\mbox{pen}(m) = \mu \Big(1+6L+8\sqrt{L}\Big)\frac{D_m}{n}=c\times\frac{D_m}{n}.
\end{equation}

In the more general context, that is in the case of irregular partitions, the numbers of models having the same dimension $D$ is exponential and 
satisfies
$$Card\Big\{ m\in \mathcal{M}, \vert m \vert=D\Big\}={n-1\choose D-1}
\leq {n\choose D}.$$
 In that case we choose $L_m$ depending on the dimension $D_m$. With $L$ depending on $D$, $\Sigma$ in \eref{sigma} satisfies \begin{eqnarray*}
\Sigma&=&\sum_{D\geq 1}e^{-L_{D}D} Card\{m\in \mathcal{M}, \vert m \vert=D\}\\
&\leq& \sum_{D\geq 1}e^{-L_{D}D} {n\choose D}\\
&\leq& \sum_{D\geq 1}e^{-L_{D}D}\Big(\frac{en}{D}\Big)^{D}\\
&\leq&\sum_{D\geq 1}e^{-D\Big(L_{D}-1-\log{(\frac{n}{D})}\Big)}
\end{eqnarray*}
So taking $L_{D}=2+\log{(\frac{n}{D})}$ leads to $\Sigma <\infty$ and the penalty becomes
\begin{equation}\label{pen}
\mbox{pen}(m)=\mu\times \mbox{pen}_{\mbox{shape}}(m),
\end{equation}
where
\begin{equation}
\mbox{pen}_{\mbox{shape}}(m)=\frac{D_m}{n}\Big[ 13+ 6\log{\Big(\frac{n}{D_m}\Big)}+8\sqrt{2+\log{\Big(\frac{n}{D_m}\Big)}}\Big].
\end{equation}
 The constant  $\mu$ can be calibrated using the slope heuristics  Birg\'e and Massart \citeyear{birge2007} (see Section \ref{slope}).

\begin{remark}
In Theorem~\ref{theo2},  we do not assume that the target function $f_0$ is piecewise constant. However in many contexts, for instance in segmentation, we might want to consider that
$f_0$ is piecewise constant or can be well approximated by piecewise constant functions. That means  there exists of partition of  $\mathcal{X}$ within which the observations follow the same distribution and between which observations have different distributions. 
\end{remark}

\section{Simulations}\label{S4}
 In this section we present numerical simulation to  study the non-asymptotic properties of the model selection procedure introduced in Section~\ref{modsel}. 
More precisely, the numerical properties of the estimators built by  model selection with our criteria  are compared with
those of the estimators resulting from model selection using the well known criteria AIC and BIC.

%
%

\subsection{Simulations frameworks}
We consider the model defined in~\eref{model} with $f_0 : [0,1] \rightarrow \mathbb{R}$. The aim is to estimate $f_0$. We consider  the collection of models $(S_m)_{m \in\mathcal{M}}$, where 
$$S_m=\mbox{Vect}\{\ind_{[\frac{k-1}{D_m},\frac{k}{D_m}[} ~\mbox{such that}~1\leq k \leq D_m\},$$
and
$\mathcal{M}$ is the collection of regular partitions $$m=\left\{\Big [\frac{k-1}{D_m},\frac{k}{D_m}\Big[, \mbox{ such that } 1\leq k\leq D_m, \right\},$$
\mbox{ where } $$ D_m \leq \frac{n}{\log n}.$$
The collection of estimators is defined in Lemma~\ref{Projhisto}.  Let us thus consider four penalties.
\begin{itemize}
\item the AIC criretion defined by
$$ \mbox{pen}_{\mbox{AIC}}=\frac{D_m}{n};$$
\item the BIC criterion defined by
$$ \mbox{pen}_{\mbox{BIC}}=\frac{\log n }{2n}D_m;$$
\item the penalty proportional to the dimension as in~\eref{penlin} defined by
$$ \mbox{pen}_{\mbox{lin}}=c\times\frac{D_m}{n};$$
\item and the penalty defined in~\eref{pen} by
$$ \mbox{pen}= \mu\times \mbox{pen}_{\mbox{shape}}(m).$$
\end{itemize}
$\mbox{pen}_{\mbox{lin}}$ and $\mbox{pen}$  are penalties depending on some unknown multiplicative constant (c and $\mu$ respectively) to be calibrated. As previously said we will use the "slope heuristics" introduced in~Birg\'ea nd Massart \citeyear{birge2007} to calibrate the multiplicative constant.
We have distinguished two cases:
\begin{itemize}
\item The case where there exists $m_o\in \mathcal{M}$ such that  the true function belong to $S_{m_{o}}$ \textit{i.e.} where $f_0$ is piecewise constant,
\begin{eqnarray}
\notag \mbox{Mod1:}~~ f_0&=&0.5\ind_{[0,1/3)}+\ind_{[1/3,0.5)}+2\ind_{[0.5,2/3)}+0.25\ind_{[2/3,1]}\\
\notag  \mbox{Mod2:}~~f_0&=&0.75\ind_{[0,1/4]} +0.5\ind_{[1/4,0.5)}+0.2\ind_{[0.5,3/4)}+0.3\ind_{[3/4,1]}.
\end{eqnarray}
\item The second case, $f_0$ does not belong to any $S_m$, $m\in \mathcal{M}$ and is chosen in the following way:
\begin{eqnarray}
\notag \mbox{Mod3:} ~~f_0(x)&=&\sin{(\pi x)} \\
\notag \mbox{Mod4:}~~ f_0(x)&=&\sqrt{x}.
\end{eqnarray}
\end{itemize}

In each case, the $x_i$'s are simulated according to uniform distribution on $[0,1].$
 
The Kullback-Leibler divergence is definitely not suitable to evaluate the quality  of an estimator. Indeed, given a model  $S_m$, there is a positive probability that on one of the interval $I\in m$ we have $\pi_{\hat f_{ m}}^{(I)}=0$ or $\pi_{\hat f_{ m}}^{(I)}=1$, which implies that $\mathcal{K}(\pi_{f_0}^{(n)},\pi_{\hat f_{m}}^{(n)})=+\infty$. So we will use the Hellinger distance to evaluate the quality of an estimator. 

Even if an oracle inequality seems of  no practical use, it  can serve as a benchmark to evaluate the performance of any data driven selection procedure.
Thus  model selection performance of each procedure is evaluated by the following benchmark
\begin{equation}
C^{*}:=\frac{\mathbb{E}\Big[h^{2}(\mathbb{P}_{f_0}^{(n)},\mathbb{P}_{\hat f_{\hat m}}^{(n)})\Big]}{\mathbb{E}\Big[\inf_{m\in \mathcal{M}}h^{2}(\mathbb{P}_{f_0}^{(n)},\mathbb{P}_{\hat f_{m}}^{(n)})\Big]}.
\end{equation}
$C^{*}$ evaluate how far is the selected estimator to the oracle. 
The values of $C^{*}$ evaluated for each procedure with different  sample size $n\in \{100, 200,\dots,1000\}$ 
are reported in Figure~\ref{fig1} , Figure~\ref{fig2}, Figure~\ref{fig1_1} and Figure~\ref{fig2_1}. For each sample size $n\in \{100, 200,\dots,1000\}$, the expectation was estimated using mean over 1000  simulated datasets. 
%
  \subsection{Slope heuristics }\label{slope}
The aim of this section is to show how the penalty in Theorem~\ref{theo2} can be calibrated in practice using the main ideas of data-driven penalized model selection criterion proposed by Birg\'e and Massart \citeyear{birge2007}. 
We calibrate penalty using "slope heuristics"  first introduced and theoretically validated by~Birg\'e and Massart \citeyear{birge2007} in a gaussian homoscedastic setting.  Recently it has also been theoretically  validated  in the heteroscedastic random-design case by Arlot \citeyear{arlot2009}  and for least squares density estimation by~Lerasle \citeyear{lerasle}. Several encouraging applications of this method are developed in many other frameworks  (see   for instance in clustering and variable selection for categorical multivariate  data ~Bontemps and Toussile \citeyear{bontemps}, for variable selection and clustering via Gaussian mixtures Maugis and Michel \citeyear{maugis2011},  in multiple change points detection Lebarbier ~\citeyear{lebarbier}). Some overview and implementation of the slope heuristics can be find in Baudry \textit{et al.}~\citeyear{baudry}. 

We now describe the main idea of those heuristics, starting from that main goal of the model selection, that is to choose the best estimator of $f_0$ among a collection of estimators $\{\hat f_m\}_{m \in \mathcal{M}}$. Moreover, we expect that this best estimator mimics the so-called oracle defined as (\ref{m*}). 
To this aim, the great challenge is to build a penalty  function such that the selected model $\hat m$ is 
nearly as good as the oracle. In the following we call the ideal penalty the penalty that leads to the choice of $m*$.
Using that $$\mathcal{K}(\mathbb{P}_{f_0}^{(n)},\mathbb{P}_{\hat f_m}^{(n)})=\gamma(\hat f_m)-\gamma(f_0),$$ then, by definition, $m*$
defined in \eref{m*} satisfies
$$m*=\arg\min_{m\in \mathcal{M}}[\gamma(\hat f_m)-\gamma(f_0)]=\arg\min_{m\in \mathcal{M}}\gamma(\hat f_m).$$
  The ideal penalty, leading to the choice of the oracle $m*$, is thus
 $[\gamma(\hat f_m)-\gamma_n(\hat f_m)]$,  for $
m\in \mathcal{M}.$ As the matter of fact, by replacing $\mbox{pen}_{id}(\hat f_m)$ by its value, we obtain
\begin{eqnarray*}
\arg\min_{m \in \mathcal{M}}[  \gamma_n(\hat f_m)+\mbox{pen}_{id}(\hat f_m)]&=&\arg\min_{m \in \mathcal{M}}[  \gamma_n(\hat f_m)+
\gamma(\hat f_m)-\gamma_n(\hat f_m)]\\&=&\arg\min_{m\in \mathcal{M}}[\gamma(\hat f_m)]\\
&=&m*.\end{eqnarray*}
Of course this ideal  penalty always selects the oracle model but depends on the unknown function $f_0$ throught the sample distribution, since
$\gamma(t)=\mathbb{E}_{f_0}[\gamma_n(t)].$ A natural idea is to choose $\mbox{pen}(m)$ as close as possible  to $\mbox{pen}_{id}(m)$ for every $m\in \mathcal{M}$. Now, we use that this ideal penalty can be decomposed into
\begin{eqnarray*}
\mbox{pen}_{id}(m)&=&\gamma(\hat f_m) -\gamma_n(\hat f_m) 
= v_m+\hat v_m+ e_m,
\end{eqnarray*} 
where $$v_m= \gamma(\hat f_m)-\gamma(f_m), \quad \hat v_m=\gamma_n(f_m)-\gamma_n(\hat f_m), \mbox{ and } \quad e_m=\gamma(f_m)-\gamma_n(f_m).$$

The slope heuristics relies on two points: 
\begin{itemize}
\item The existence of a minimal penalty $\mbox{pen}_{\mbox{min}}(m)=\hat v_m$ such that  when the penalty is smaller than $\mbox{pen}_{\mbox{min}}$ the selected model is one of the most complex models. Whereas, penalties  larger than  $\mbox{pen}_{\mbox{min}}$  lead to a selection of   models with "reasonable" complexity.
\item  Using concentration arguments, it is reasonable to consider that uniformly over $\mathcal{M}$, $\gamma_n(f_m)$ is close to its expectation 
which implies that $e_m\approx 0$.   In the same way, since $\hat v_m$ is a empirical version of $v_m$, it is also reasonable to consider that  $v_m\approx \hat v_m$.  Ideal penalty is thus approximately given by $2 \hat v_m$, and thus
\begin{eqnarray*}
\mbox{pen}_{id}(m)
&\approx &2 \mbox{pen}_{min}(m).
\end{eqnarray*}
\end{itemize} 
In practice, $\hat v_m$ can be estimated from the data provided that ideal penalty $\mbox{pen}_{id}(.)=\kappa_{id}\mbox{pen}_{shape}(.)$ is known up to a multiplicative factor. A major point of the slope heuristics is that
$$\frac{\kappa_{id}}{2}\mbox{pen}_{shape}(.)$$
is a good estimator of  $\hat v_m$ and this provides the minimal penalty. 

 Provided that $\mbox{pen}=\kappa\times \mbox{pen}_{shape}$ is known up to a multiplicative constant $\kappa$  that is to be calibrated,  we combine the previously heuristic to the method usually known  as dimension jump method. In practice, we consider a grid $\kappa_1,\dots,\kappa_M$,  where each $\kappa_j$ leads to a selected model $\hat m_{\kappa_i}$ with dimension $D_{\hat m_{\kappa_i}}$. The constant $\kappa_{min}$ which corresponds to the value such that $\mbox{pen}_{min}=\kappa_{min}\times \mbox{pen}_{shape}$, is estimated using the first point of the "slope heuristics".  If $D_{\hat m_{\kappa_j}}$ is plotted as a function of $\kappa_j$, $\kappa_{min}$ is such that $D_{\hat m_{\kappa_j}}$ is "huge" for $\kappa< \kappa_{min}$ and "reasonably small" for $\kappa>\kappa_{min}$. So $\kappa_{min}$ is the value at the position of the biggest jump. For more details about this method we refer the reader to Baudry \textit{et al.}~\citeyear{baudry} and Arlot and Massart \citeyear{arlot2009}.

Figures~\ref{fig1}  and \ref{fig1_1} are the cases where the true function is  piecewise constant. Figure~\ref{fig2} and Figure~\ref{fig2_1}
are situations where the true function does not belong to any model in the given collection.  The performance of criteria depends on the sample size $n$. In  these two situations we observe that our two model selection procedures are comparable, and  their performance increases with $n$.  While the performance of model selected by BIC decreases with $n$. Our criteria outperformed the AIC for all $n$. The BIC criterion is better than our criteria for $n\leq 200$. For  $200< n\leq 400$,  the performance of the model selected by BIC is quite the same as the performance of models selected by our criteria.
Finally for $n>400$ our criteria outperformed the BIC.

  Theoretical results and simulations raise the following question :  why our criteria are better than BIC for quite large values of $n$ yet  theoretical results are non asymptotic?  To answer this question we can say that, in simulations, to calibrate our penalties we have used "slope heuristics",  and  those heuristic are based on asymptotic arguments (see Section~\ref{slope}).

%
 
\begin{figure}
\includegraphics[height=20cm,width=15cm]{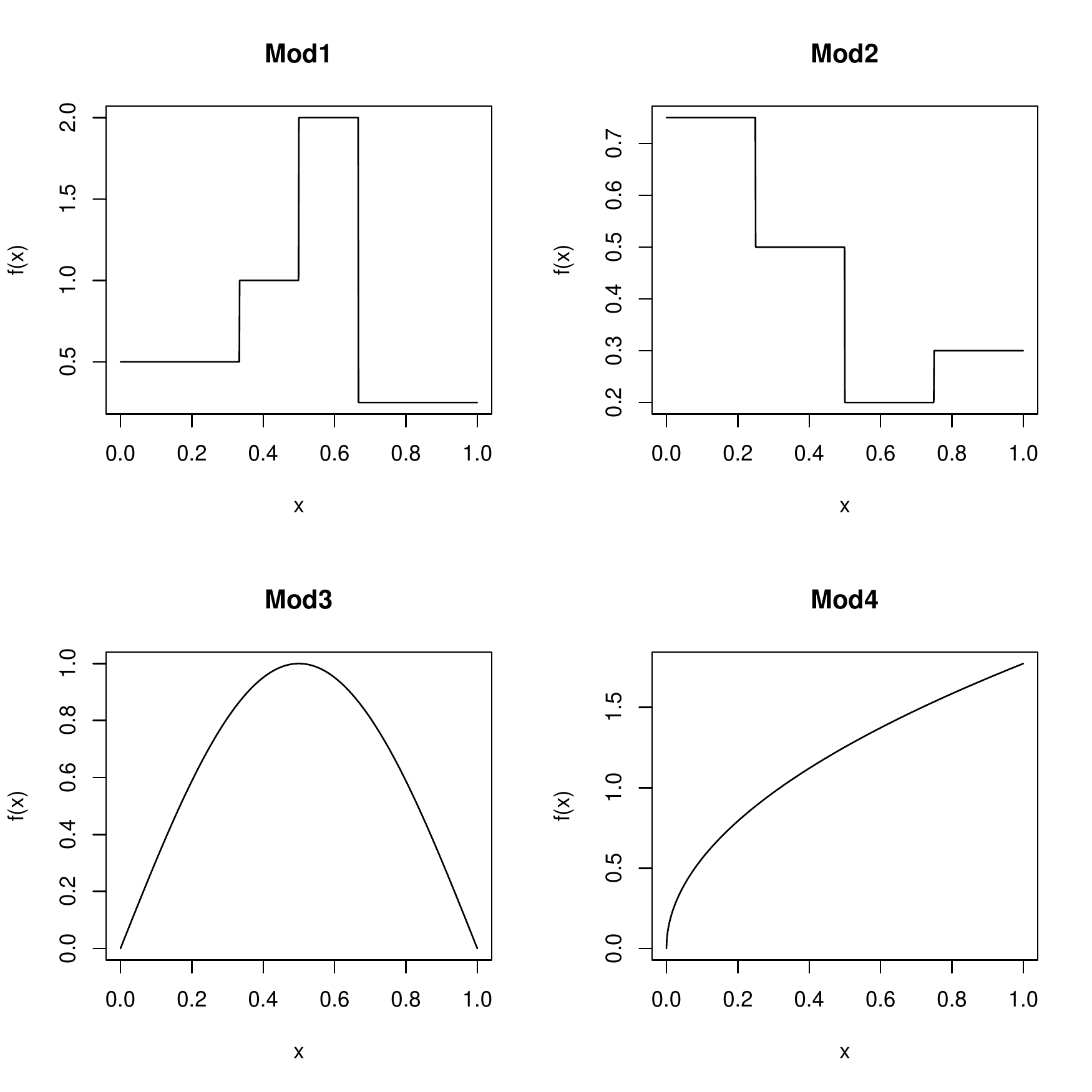}
\caption{Different functions $f_0$ to be estimated}
\end{figure}

\begin{figure}
\includegraphics[height=11cm,width=13cm]{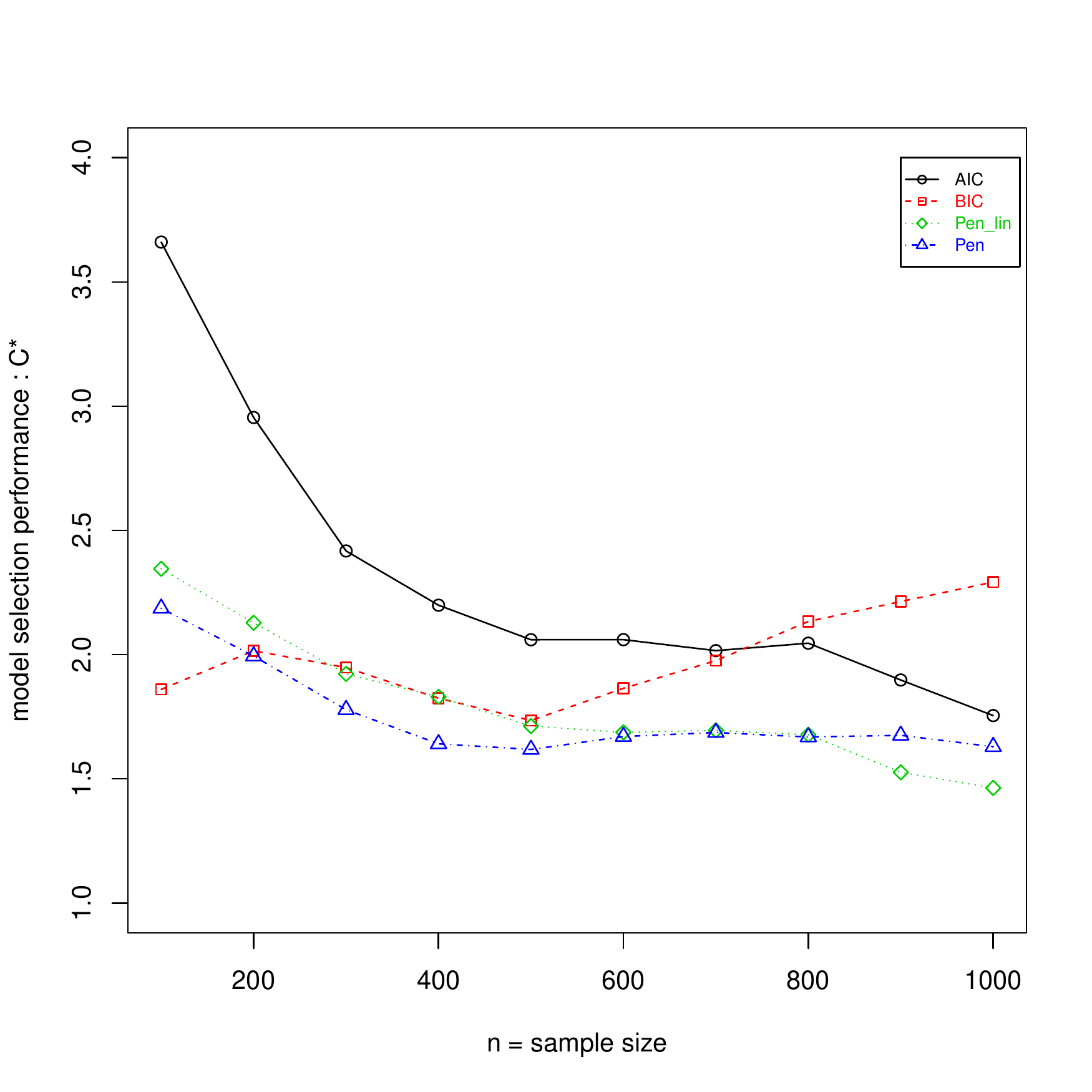}
\caption{\label{fig1} Model selection performance ($C^{*}$) as a function  of sample size n, with each penalty, Mod1.}
\end{figure}
\begin{figure}
\includegraphics[height=11cm,width=13cm]{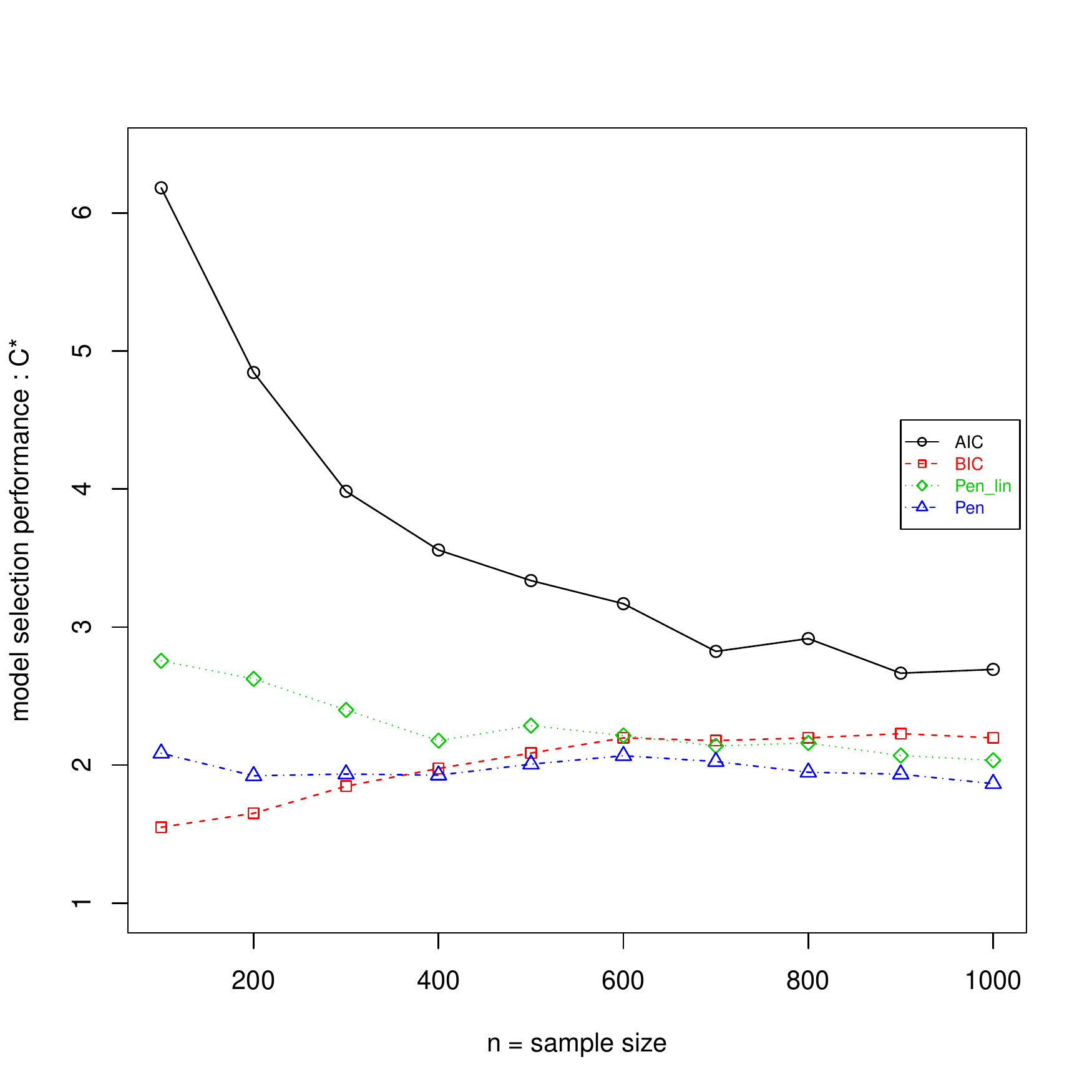}
\caption{\label{fig1_1} Model selection performance ($C^{*}$) as a function  of sample size n, with each penalty, Mod2.}
\end{figure}

\begin{figure}
\includegraphics[height=11cm,width=13cm]{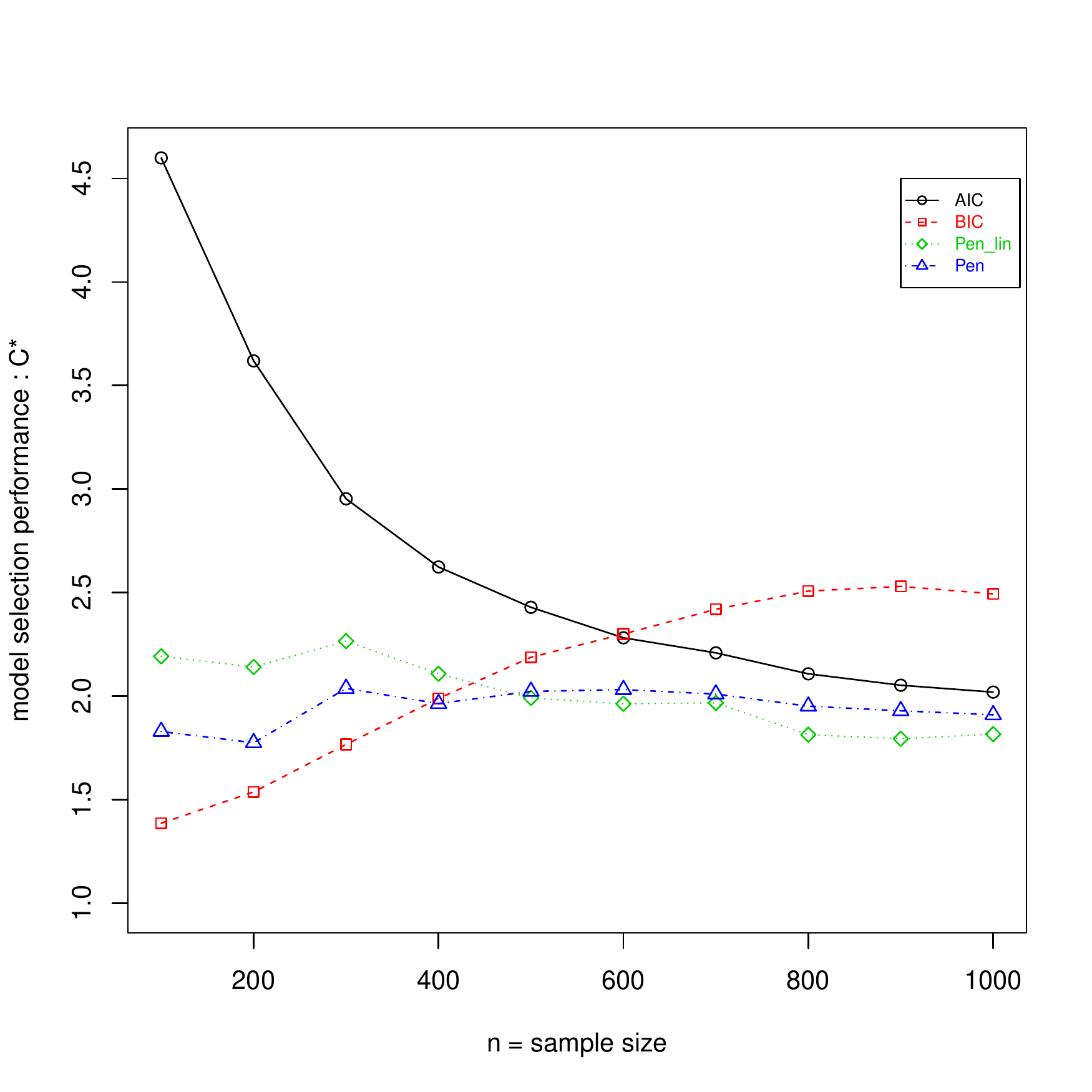}
\caption{\label{fig2} Model selection performance  ($C^{*}$) as a function   of sample size n, with each penalty, Mod3.}
\end{figure}
\begin{figure}
\includegraphics[height=11cm,width=13cm]{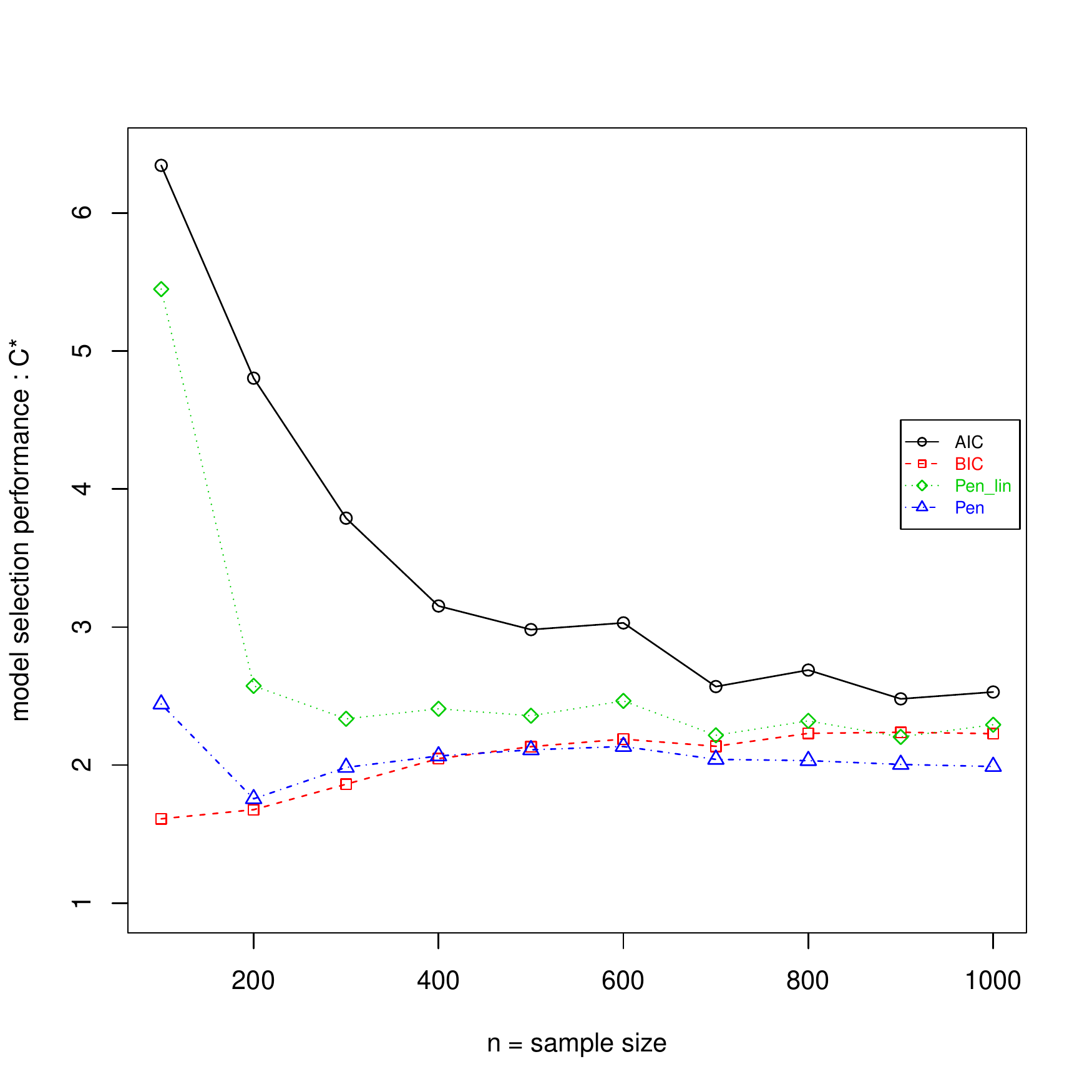}
\caption{\label{fig2_1} Model selection performance  ($C^{*}$) as a function   of sample size n, with each penalty, Mod4.}
\end{figure}


\section{Proofs}\label{S5}
\setcounter{equation}{0}
\setcounter{lem}{0}
\setcounter{theo}{0}


\subsection{Notations and technical tools}
Subsequently we will use the following notations.
Denote by
$\parallel f\parallel_n$   and $\langle f,g \rangle_n $ the empirical euclidian norm
and the  inner product 
$$\parallel f\parallel_n^2=\frac{1}{n}\sum_{i=1}^n f^2(x_i), \mbox{ and }\langle f,g\rangle_n=\frac{1}{n}\sum_{i=1}^n f(x_i)g(x_i).$$
Note that $\parallel . \parallel_n$ is a semi norm on the space $\mathcal{F}$ of functions $g : \mathcal{X}  \longrightarrow \mathbb{R}$, but is a norm in the quotient space  $\mathcal{F} \slash \mathcal{R} $ associated to the equivalence relation $\mathcal{R}$ : $g ~\mathcal{R}~ h$ if and only if $g(x_i)=h(x_i)$ for all $i\in \{1,\dots,n\}$.
It follows  from \eref{gamma_n} that  $\gamma$ defined in \eref{gamma} can be expressed as the sum of a centered empirical process and of the estimation criterion $\gamma_n$. More precisely, denoting
by $\vec{\varepsilon}=(\varepsilon_1,\cdots,\varepsilon_n)^T$, with $\varepsilon_i=Y_i-\mathbb{E}_{f_0}(Y_i),$ for all $f$, we have
\begin{eqnarray}\label{decompgamma}
\gamma(f)=\gamma_n(f)+\frac{1}{n}\sum_{i=1}^n \varepsilon_i f(x_i):= \gamma_n(f)+\langle\vec{\varepsilon},f\rangle_n.
\end{eqnarray}
 Easy calculations show that  for $\gamma$ defined in \eref{gamma} we have,
\begin{eqnarray*}
\mathcal{K}(\mathbb{P}_{f_0}^{(n)},\mathbb{P}_{f}^{(n)})&=&\frac{1}{n}\int \log\left( \frac{\mathbb{P}_{f_0}^{(n)}}{\mathbb{P}_{f}^{(n)}}\right) d\mathbb{P}_{f_0}^{(n)}
=\gamma(f)-\gamma(f_0)\\
&=&\frac{1}{n}\sum_{i=1}^n \left[ \pi_{f_0}(x_i)\log\left( \frac{\pi_{f_0}(x_i)}{ \pi_{f}(x_i) }\right) 
+(1-\pi_{f_0}(x_i))\log\left( \frac{1-\pi_{f_0}(x_i)}{ 1-\pi_{f}(x_i) }\right) 
 \right].
\end{eqnarray*}

Let us recall the  usual bounds (see Castellan \citeyear{Castellan2003}) for kullback-Leibler information:
\begin{lem}
\label{borneK}
For positive densities $p$ and $q$ with respect to $\mu$, if $f=\log(q/p)$, then
\begin{eqnarray*}
\frac{1}{2}\int f^2 (1\wedge e^f) p\, d\mu \leqslant \mathcal{K}(p,q)\leqslant \frac{1}{2}\int f^2 (1\vee e^f) p \,d\mu .
\end{eqnarray*}

\end{lem}

\subsection{Proof of Proposition \ref{borne1}:}
By definition of $\hat{f}_m$, for all $f \in S_m\cap \mathbb{L}_\infty(C_0)$, 
$\gamma_n(\hat{f}_m)-\gamma_n(f)\leqslant 0.$
We apply  \eref{decompgamma}, with $f=f_m$ and $f=\hat f_m$,
\begin{eqnarray*}
\gamma(\hat{f}_m)-\gamma(f_0)\leqslant \gamma(f_m)-\gamma(f_0)+\langle\vec{\varepsilon},\hat{f}_m-f_m\rangle_n.
\end{eqnarray*}
As usual, the main part of the proof relies on the study of the empirical process
$\langle\vec{\varepsilon},\hat{f}_m-f_m\rangle_n$.
Since $\hat{f}_m-f_m$ belongs to $S_m$,  $\hat{f}_m-f_m=\sum_{j=1}^{D_m}\alpha_j\psi_j$, where $\{\psi_1,\dots,\psi_{D_m}\},$ 
is an orthonormal basis of $S_m$ and consequently
$$\langle\vec{\varepsilon},\hat{f}_m-f_m\rangle_n= \sum_{j=1}^{D_m}\alpha_j \langle\vec{\varepsilon},\psi_j\rangle_n.$$
Applying Cauchy-Schwarz inequality we get
\begin{eqnarray*}
\langle\vec{\varepsilon},\hat{f}_m-f_m\rangle_n&\leqslant& \sqrt{\sum_{j=1}^{D_m}\alpha_j^{2}}\sqrt{\sum_{j=1}^{D_m}\left(\langle\vec{\varepsilon},\psi_j\rangle_n\right)^{2}}
\\
&=&\lVert \hat{f}_m-f_m \rVert_n\sqrt{\sum_{j=1}^{D_m}\left(\frac{1}{n}\sum_{i=1}^{n}\varepsilon_{i}\psi_j(x_i)\right)^{2}}.
\end{eqnarray*}
We now apply Lemma~\ref{lm} (See Section \ref{appendix} for the proof of Lemma~\ref{lm})

\begin{lem}\label{lm}
Let  $\mathcal{S}_m$ the model defined in \eref{model_col_gen} and $\{\psi_1,\dots,\psi_{D_m}\}$ an orthonormal basis of the linear span $\{\phi_k,~ k\in m\}$. We also denote by $\Lambda_m$  the set of $\beta=(\beta_1,...,\beta_D)$ such that  $f_\beta(.)=\sum_{j=1}^{D}\beta_j\psi_j(.)$ satisfies  $f_{\beta}\in \mathcal{S}_m\cap \mathbb{L}_\infty(C_0)$. 
Let  $\beta^{*}$ be any minimizer of the function $\beta \rightarrow \gamma(f_\beta)$ over $\Lambda_m$, we have
\begin{equation}
 \frac{\mathcal{U}_0^{2}}{2}\lVert f_{\beta}-f_{\beta^{*}}\rVert^2_{n}\leq \gamma(f_\beta)-\gamma(f_{\beta^{*}}),
\end{equation}
where $\mathcal{U}_0=e^{C_0}/(1+e^{C_0})^2$.
\end{lem}
Then we have 
$$\langle\vec{\varepsilon},\hat{f}_m-f_m\rangle_n\quad  
\leqslant \quad \sqrt{\sum_{j=1}^{D_m}\left(\langle\vec{\varepsilon},\psi_j\rangle_n\right)^{2}}\frac{\sqrt{2}}{\mathcal{U}_0} \sqrt{\gamma(\hat{f}_m)-\gamma(f_m)} $$
Now we use that for every positive numbers, $a$, $b$, $x$, $ab\leqslant (x/2)a^2+ [1/(2x)] b^2$, and infer that
\begin{equation*}
 \gamma(\hat{f}_m)-\gamma(f_0)\leq
\gamma(f_m)-\gamma(f_0)+ \frac{x}{\mathcal{U}_0^{2}}\sum_{j=1}^{D_m}\left(\langle\vec{\varepsilon},\psi_j\rangle_n\right)^{2}  +
(1/2x)( \gamma(\hat{f}_m)-\gamma(f_m)).
\end{equation*}
For $x>1/2$, it follows that
\begin{eqnarray*}
\mathbb{E}_{f_0}[\gamma(
\hat{f}_m)-\gamma(f_0)]\leqslant  
\gamma(f_m)-\gamma(f_0)+\frac{2x^{2}}{(2x-1)\mathcal{U}_0^{2}} \mathbb{E}_{f_0}\left[\sum_{j=1}^{D_m}\left(\langle\vec{\varepsilon},\psi_j\rangle_n\right)^{2}\right] .
\end{eqnarray*}
 We conclude the proof by using  that
$$\mathbb{E}_{f_0}\left[\sum_{j=1}^{D_m}\left(\langle\vec{\varepsilon},\psi_j\rangle_n\right)^{2}\right]\leqslant \frac{D_m}{4n}. $$
\qed

\subsection{Proof of Theorem~\ref{theo1}}
By definition, for all $m \in \mathcal{M}$,
$$\gamma_n(\hat f_{\hat{m}})+\mbox{pen}(\hat{m})\leqslant \gamma_n(\hat{f_m})+\mbox{pen}(m) \leqslant \gamma_n(f_m)+\mbox{pen}(m).$$
Applying  \eref{decompgamma} we have
\begin{equation}\label{decomp1}
\mathcal{K}(\mathbb{P}_{f_0}^{(n)},\mathbb{P}_{\hat{f}_{\hat m}}^{(n)})\leqslant \mathcal{K}(\mathbb{P}_{f_0}^{(n)},\mathbb{P}_{f_m}^{(n)})+\langle \vec{\varepsilon},\hat{f}_{\hat m}- f_m\rangle_n+\mbox{pen}(m)-\mbox{pen}(\hat{m}).
\end{equation}
It remains to study $\langle \vec{\varepsilon}, \hat{f}_{\hat m}-f_m\rangle_n$, using the following lemma, which is a modification of Lemma 1 in Durot \textit{et al.} \citeyear{DurotLebarbierTocquet}.
\begin{lem}\label{control}
 For every $D$, $D^{\prime}$ and $x\geqslant 0$ we have
$$\mathbb{P}\left(\sup_{u\in \big(S_{D}\cap \mathbb{L}_\infty(C_0)+S_{D^{\prime}}\cap \mathbb{L}_\infty(C_0)\Big)}\frac{\langle\vec{\varepsilon} ,u\rangle_n}{\parallel u \parallel_n}\rangle\sqrt{\frac{ D+D^{\prime}}{4n}}+\sqrt{\frac{5x}{n}}\right)\leqslant \exp{(-x)}.$$
Fix $\xi>0$ and let $\Omega_\xi(m)$ denote the  event
$$\Omega_\xi(m)=\bigcap_{m^{\prime}\in \mathcal{M}}\left\{ \sup_{u\in \Big(S_{m}\cap \mathbb{L}_\infty(C_0)+S_{m^{\prime}}\cap \mathbb{L}_\infty(C_0)\Big)}\frac{\langle\vec{\varepsilon} ,u\rangle_n}{\parallel u \parallel_n} \leq
\sqrt{\frac{ D_m+D_{m^{\prime}}}{4n}}+\sqrt{5(L_{m^{\prime}}D_{m^{\prime}}+\xi)/n} \right\}.$$
Then we have
\begin{equation}\label{prob_sigma}
\mathbb{P}\left(\Omega_\xi(m)\right)\geqslant 1-\Sigma\exp(-\xi).
\end{equation}
\end{lem}
See  the Appendix for the proof of this lemma.
Fix $\xi>0$, applying Lemma \ref{control}, we infer that on the event $\Omega_\xi(m),$
\begin{align*}
 \langle \vec{\varepsilon}, \hat{f}_{\hat m}-f_m\rangle_n&
\leqslant \left(\sqrt{\frac{ D_m+D_{\hat{m}}}{4n}}+\sqrt{5\frac{L_{\hat{m}}D_{\hat{m}}+\xi}{n}}\right)\parallel \hat{f}_{\hat m}-f_m\parallel_{n}\\ 
&\leqslant \left(\sqrt{\frac{ D_m+D_{\hat{m}}}{4n}}+\sqrt{5\frac{L_{\hat{m}}D_{\hat{m}}+\xi}{n}}\right)\left(\parallel \hat{f}_{\hat m}-f_{0}\parallel_n+\parallel f_{0}-f_m\parallel_{n}\right)\\
&\leqslant \left(\sqrt{D_{\hat{m}}}\left(\frac{1}{\sqrt{4n}}+\sqrt{\frac{5L_{\hat{m}}}{n}}\right)+\sqrt{\frac{D_m}{4n}}+\sqrt{5\frac{\xi}{n}}\right)\left(\parallel \hat{f}_{\hat m}-f_{0}\parallel_n+\parallel f_{0}-f_m\parallel_{n}\right).
\end{align*}
Applying  that $2xy\leqslant\theta x^{2}+ \theta^{-1}y^{2}$, for all $ x>0$, $y>0$, $\theta>0$, we get that on $\Omega_\xi(m)$ and for every $\eta\in ]0,1[$
\begin{eqnarray*}
  \langle \vec{\varepsilon}, \hat{f}_{\hat m}-f_m\rangle_n\!\!\!&\leqslant &\!\!\!(\frac{1-\eta}{2})\left[(1+\eta)\parallel \hat{f}_{\hat m}-f_{0}\parallel_{n}^{2}+(1+\eta^{-1})\parallel f_{0}-f_m\parallel_{n}^{2}\right]\\
&+&\frac{1}{2(1-\eta)}\left[(1+\eta)
D_{\hat{m}}\left(\frac{1}{\sqrt{4n}}+\sqrt{\frac{5L_{\hat{m}}}{n}}\right)^{2}+(1+\eta^{-1})\left(\sqrt{\frac{D_m}{4n}}+\sqrt{\frac{5\xi}{n}}\right)^{2}\right]\\
\!\!\!&\leqslant&\!\!\! \frac{1-\eta^{2}}{2}\parallel \hat{f}_{\hat m}-f_{0}\parallel_{n}^{2}+
\frac{\eta^{-1}-\eta}{2}\parallel f_{0}-f_m\parallel_{n}^{2}+\frac{1+\eta}{2(1-\eta)}D_{\hat{m}}\left(\frac{1}{\sqrt{4n}}+\sqrt{\frac{5L_{\hat{m}}}{n}}\right)^{2}\\
&&+\frac{1+\eta^{-1}}{1-\eta}\Big(\frac{D_m}{4n}+\frac{5\xi}{n}\Big).
\end{eqnarray*}
If $\mbox{pen}(m)\geqslant \Big(\lambda D_m\left(\frac{1}{2}+\sqrt{5L_m}\right)^{2}\Big)/n, $ with $\lambda>0$, we have 
\begin{eqnarray*}
  \langle \vec{\varepsilon}, \hat{f}_{\hat m}-f_m\rangle_n
\!\!\!&\leqslant&\!\!\! \frac{1-\eta^{2}}{2}\parallel \hat{f}_{\hat m}-f_{0}\parallel_{n}^{2}+
\frac{\eta^{-1}-\eta}{2}\parallel f_{0}-f_m\parallel_{n}^{2}+\frac{1+\eta}{2(1-\eta)\lambda}\mbox{pen}(\hat{m})+\frac{1+\eta^{-1}}{(1-\eta)\lambda}\mbox{pen(m)}\\
&&+\frac{1+\eta^{-1}}{1-\eta}\frac{5\xi}{n}.
\end{eqnarray*}

It follows from (\ref{decomp1}) that
\begin{eqnarray*}
 \mathcal{K}(\mathbb{P}_{f_0}^{(n)}, \mathbb{P}_{\hat{f}_{\hat m}}^{(n)})&\leqslant&\mathcal{K}(\mathbb{P}_{f_0}^{(n)},\mathbb{P}_{f_m}^{(n)})+\frac{1-\eta^{2}}{2}\parallel \hat{f}_{\hat m}-f_{0}\parallel_{n}^{2}+\frac{\eta^{-1}-\eta}{2}\parallel f_{0}-f_m\parallel_{n}^{2}\\&&+\frac{1+\eta}{2(1-\eta)\lambda}\mbox{pen}(\hat{m})+\frac{1+\eta^{-1}}{(1-\eta)\lambda}\mbox{pen(m)}+\frac{1+\eta^{-1}}{1-\eta}\frac{5\xi}{n}+\mbox{pen}(m)-\mbox{pen}(\hat{m}).
\end{eqnarray*}
Taking
$\lambda=(\eta+1)/(2(1-\eta))$,
we have
\begin{multline*}
\mathcal{K}(\mathbb{P}_{f_0}^{(n)},\mathbb{P}_{\hat{f}_{\hat m}}^{(n)})\leqslant \mathcal{K}(\mathbb{P}_{f_0}^{(n)},\mathbb{P}_{f_m}^{(n)})\\+\frac{4\lambda}{(2\lambda+1)^{2}}\parallel \hat{f}_{\hat m}-f_{0}\parallel_{n}^{2}+
  \frac{4\lambda}{4\lambda^{2}-1}\parallel f_{0}-f_m\parallel_{n}^{2}+\frac{6\lambda+1}{2\lambda-1}\mbox{pen}(m)+\frac{10\lambda(2\lambda+1)}{2\lambda-1}\frac{\xi}{n}.
\end{multline*}
Now we use the following lemma  (see Lemma 6.1 in Kwemou \citeyear{kwemou}) that allows to connect empirical norm and Kullback-Leibler divergence.
\begin{lem}\label{l81}
 Under Assumptions~(\ref{A1}), for all $m\in \mathcal{M}$ and all $t\in S_m\cap \mathbb{L}_\infty(C_0)$, we have
$$c_{min}\lVert t-f_{0}\rVert_{n}^{2}\leqslant \mathcal{K}(\mathbb{P}_{f_0}^{(n)},\mathbb{P}_{t}^{(n)}) \leqslant c_{max}\lVert t-f_{0}\rVert_{n}^{2}.$$
where $c_{min}$ and $c_{max}$ are constants depending on $C_{0}$ and $c_{1}.$
\end{lem}
Consequently
$$\mathcal{K}(\mathbb{P}_{f_0}^{(n)},\mathbb{P}_{\hat{f}_{\hat m}}^{(n)})\leqslant C(c_{min})\left\{\mathcal{K}(\mathbb{P}_{f_0}^{(n)},\mathbb{P}_{f_m}^{(n)}) +\mbox{pen}(m)\right\}+C_1(c_{min})\frac{\xi}{n},$$
where 
$$C(c_{min})=\max\left\{ \frac{1+\frac{4\lambda}{(4\lambda^{2}-1)c_{min}}}{1-\frac{4\lambda}{c_{min}(2\lambda+1)^{2}}};
\frac{\frac{6\lambda+1}{2\lambda-1}}{1-\frac{4\lambda}{c_{min}(2\lambda+1)^{2}}}\right\}~\mbox{and}
~~C_1(c_{min})=\frac{\frac{10\lambda(2\lambda+1)}{2\lambda-1}}{1-\frac{4\lambda}{c_{min}(2\lambda+1)^{2}}}.$$
Thus we take $\lambda$  such that 
\begin{equation}\label{cond_lambda}
1-\frac{4\lambda}{c_{min}(2\lambda+1)^{2}}>0,
\end{equation}
where $c_{min}$ depends on the bound of  the true function $f_0$.
By definition of $\Omega_\xi(m)$ and \eref{prob_sigma}, there exists a random variable $V\geqslant 0$ with $\mathbb{P} (V>\xi)\leqslant \Sigma\exp{(-\xi)}$ and  $\mathbb{E}_{f_0}(V)\leqslant \Sigma,$
such that
$$\mathcal{K}(\mathbb{P}_{f_0}^{(n)},\mathbb{P}_{\hat{f}_{\hat m}}^{(n)})\leqslant C(c_{min})\left\{\mathcal{K}(\mathbb{P}_{f_0}^{(n)},\mathbb{P}_{f_m}^{(n)})+\mbox{pen}(m)\right\}+C_1(c_{min})\frac{V}{n},$$
which implies that for all $m\in \mathcal{M}$,

$$\mathbb{E}_{f_0}[\mathcal{K}(\mathbb{P}_{f_0}^{(n)},\mathbb{P}_{\hat{f}_{\hat m}}^{(n)})]\leqslant C(c_{min})\left\{\mathcal{K}(\mathbb{P}_{f_0}^{(n)},\mathbb{P}_{f_m}^{(n)})+\mbox{pen}(m)\right\}+C_1(c_{min})\frac{\Sigma}{n}.$$
This concludes the proof. \qed

\subsection{Proof of Proposition \ref{borne2}:}
{Let $f_m$, $\hat f_m$, $\pi_{f_m}$ and $ \pi_{\hat f_m}$ given in Lemma~\ref{Projhisto}, proved in appendix. 
In the following, $D_m=\vert m\vert.$
For $\delta>0$, let $\Omega_m(\delta)$ be the event
\begin{eqnarray}
\label{Omega}
\Omega_m(\delta)
=
\bigcap_{J\in m} \left\lbrace \left\vert \frac{\pi_{\hat f_m}^{(J)}}{ \pi_{ f_m}^{(J)} } -1\right\vert \leqslant \delta\right\rbrace    \bigcap \left\lbrace  \left\vert \frac{1-\pi_{\hat f_m}^{(J)}}{ 1-\pi_{ f_m}^{(J)} } -1\right\vert 
\leqslant \delta \right\rbrace.\end{eqnarray}
According to pythagore's type identity and Lemma \ref{Projhisto} we write 
\begin{eqnarray*}
\mathcal{K}(\mathbb{P}_{f_0}^{(n)},\mathbb{P}_{\hat{f_{m}}}^{(n)})=\mathcal{K}(\mathbb{P}_{f_0}^{(n)},\mathbb{P}_{f_m}^{(n)})+\mathcal{K}(\mathbb{P}_{f_m}^{(n)},\mathbb{P}_{\hat{f_{m}}}^{(n)})\ind_{\Omega_m(\delta)}+
\mathcal{K}(\mathbb{P}_{f_m}^{(n)},\mathbb{P}_{\hat{f_{m}}}^{(n)})\ind_{\Omega_m^c(\delta)},
\end{eqnarray*}
where
\begin{eqnarray}\label{k_m_m}
 \mathcal{K}(\mathbb{P}_{f_m}^{(n)},\mathbb{P}_{\hat{f_{m}}}^{(n)})\!\!\!&=&\!\!\!\frac{1}{n}\sum_{i=1}^n \left[\pi_{f_m}(x_i)\log\left(\frac{\pi_{f_m}(x_i)}{\pi_{\hat f_m}(x_i)}\right)
+(1-\pi_{f_m}(x_i))\log\left(\frac{1-\pi_{f_m}(x_i)}{1-\pi_{\hat f_m}(x_i)}\right)\right]\\
\notag \!\!\!&=&\!\!\!
\frac{1}{n}\sum_{J\in m}  \vert J \vert \left[\pi_{f_m}^{(J)}\log\left(\frac{\pi_{f_m}^{(J)}}{\pi_{\hat f_m}^{(J)}}\right)
+(1-\pi_{f_m}^{(J)})\log\left(\frac{1-\pi_{f_m}^{(J)}}{1-\pi_{\hat f_m}^{(J)}}\right)\right].
\end{eqnarray}
The first step consists in showing that 
\begin{eqnarray}
\label{step1}
\frac{1-\delta}{2(1+\delta)^2}\mathcal{X}_m^2 \ind_{\Omega_m(\delta)}\leqslant \mathcal{K}(\mathbb{P}_{f_m}^{(n)},\mathbb{P}_{\hat{f_{m}}}^{(n)}) \ind_{\Omega_m(\delta)} 
\leqslant \frac{1+\delta}{2(1-\delta)^2}\mathcal{X}_m^2\ind_{\Omega_m(\delta)} ,
\end{eqnarray}
where
\begin{eqnarray}\label{chideux}
\mathcal{X}_m^2=\frac{1}{n}\sum_{J\in m} \frac{(\sum_{k\in J}\varepsilon_k)^2}{\vert J\vert\pi_{f_m}^{(J)}[1-\pi_{f_m}^{(J)}]},
\mbox{ with }\qquad 
\frac{4\rho^2 D_m}{n} \leqslant \mathbb{E}_{f_0}[ \mathcal{X}_m^2]\leqslant \frac{2D_m}{n}.
\end{eqnarray}
The second step relies on the proof of
\begin{eqnarray}
\label{step2}
\big\vert \mathbb{E}_{f_0}\left(\mathcal{K}(\mathbb{P}_{f_m}^{(n)},\mathbb{P}_{\hat{f_{m}}}^{(n)})\ind_{\Omega_m^c(\delta)}\right)\Big\vert \leqslant 2\log\left(    \frac{1}{\rho}\right)\mathbb{P}[ \Omega_m^c(\delta)] .
\end{eqnarray}
The last step  consists in showing  that for $\epsilon >0$,  since  for all $J\in m$, $\vert J\vert \geq \Gamma [\log(n)]^2$, where $\Gamma>0$ is an absolute constant, then we have
\begin{eqnarray}
\label{step3}\mathbb{P}[ \Omega_m^c(\delta)] \leqslant 4\vert m\vert \exp\left(-\frac{\delta^2}{2(1+\delta/3)}\rho^2 \Gamma [\log(n)]^2 \right)\leq\frac{\kappa(\rho,\delta,\Gamma,\epsilon)}{n^{(1+\epsilon)}}.
\end{eqnarray}  
Gathering \eref{step1}-\eref{step3}, we conclude that
\begin{eqnarray*}
\mathbb{E}_{f_0}[\mathcal{K}(\mathbb{P}_{f_0}^{(n)},\mathbb{P}_{\hat{f_{m}}}^{(n)})]&\leqslant& \mathcal{K}(\mathbb{P}_{f_0}^{(n)},\mathbb{P}_{f_{m}}^{(n)})+\frac{(1+\delta)\vert m\vert}{(1-\delta)^2n} +2\log\left(    \frac{1}{\rho}\right)\mathbb{P}[ \Omega_m^c(\delta)]\\
&\leqslant& \mathcal{K}(\mathbb{P}_{f_0}^{(n)},\mathbb{P}_{f_{m}}^{(n)})+\frac{(1+\delta)\vert m\vert}{(1-\delta)^2n} +\frac{\kappa(\rho,\delta,\Gamma,\epsilon)}{n^{(1+\epsilon)}}.
\end{eqnarray*}
We finish by proving  \eref{step1}, \eref{chideux}, \eref{step2} and \eref{step3}.
}
\paragraph{$\bullet$ Proof of \eref{step1} and \eref{chideux} :}
{
Arguing as in Castellan \citeyear{Castellan2003} and using Lemma \ref{borneK} we have
\begin{eqnarray*}
 \mathcal{K}(\mathbb{P}_{f_m}^{(n)},\mathbb{P}_{\hat{f_{m}}}^{(n)}) \!\!\!&\geqslant&\!\!\!\frac{1}{2n} \sum_{J\in m}\vert J\vert\left[\pi_{f_m}^{(J)}\left(1\wedge \frac{\pi_{\hat f_m }^{(J)} }{\pi_{f_m} ^{(J)} }   \right)\log^2\left(\frac{\pi_{f_m}^{(J)} }{\pi_{\hat f_m}^{(J)} }\right)
+(1-\pi_{f_m}^{(J)} ) \left(1\wedge \frac{1-\pi_{\hat f_m }^{(J)} }{1-\pi_{f_m} ^{(J)} }   \right)\log^2\left(\frac{1-\pi_{f_m}^{(J)}}{1-\pi_{\hat f_m}^{(J)} }\right)\right] \\
\end{eqnarray*}
and
\begin{eqnarray*}
 \mathcal{K}(\mathbb{P}_{f_m}^{(n)},\mathbb{P}_{\hat{f_{m}}}^{(n)}) \!\!\!&\leqslant& \!\!\!
\frac{1}{2n} \sum_{J\in m} \vert J\vert\left[\pi_{f_m}^{(J)} \left(1\vee \frac{\pi_{\hat f_m }^{(J)} }{\pi_{f_m}^{(J)} }   \right)\log^2\left(\frac{\pi_{f_m}^{(J)} }{\pi_{\hat f_m}^{(J)} }\right)
+(1-\pi_{f_m}^{(J)} ) \left(1\vee \frac{1-\pi_{\hat f_m }^{(J)} }{1-\pi_{f_m} ^{(J)} }   \right)\log^2\left(\frac{1-\pi_{f_m}^{(J)} }{1-\pi_{\hat f_m}^{(J)} }\right)\right]. \\
\end{eqnarray*}
It follows that
\begin{eqnarray}\label{encadrK}
 \frac{1-\delta}{2} V^2( \pi_{f_m},\pi_{\hat f_m}) \ind_{\Omega_m(\delta)}\leqslant \mathcal{K}(\mathbb{P}_{f_m}^{(n)},\mathbb{P}_{\hat{f_{m}}}^{(n)}) \ind_{\Omega_m(\delta)}\leqslant \frac{1+\delta}{2} V^2( \pi_{f_m},\pi_{\hat f_m})\ind_{\Omega_m(\delta)},
\end{eqnarray}
where $V^2( \pi_{f_m},\pi_{\hat f_m})$ is defined by
\begin{multline}\label{V}
 V^2( \pi_{f_m},\pi_{\hat f_m})=\frac{1}{n} \sum_{J\in m} \vert J\vert \frac{[\pi_{\hat f_m}^{(J)} -\pi_{f_m}^{(J)} ] ^2}{\pi_{f_m}^{(J)} }\left[ \frac{\log[\pi_{\hat f_m}^{(J)}/\pi_{f_m}^{(J)}] }{\pi_{\hat f_m}^{(J)}/\pi_{f_m}^{(J)} -1    } \right]^2\\
+\frac{1}{n} \sum_{J\in m} \vert J\vert 
\frac{[\pi_{\hat f_m}^{(J)} -\pi_{f_m}^{(J)} ] ^2}{1-\pi_{f_m}^{(J)} }\left[ \frac{\log[(1-\pi_{\hat f_m}^{(J)})/(1-\pi_{f_m}^{(J)})] }{(1-\pi_{\hat f_m}^{(J)})/(1-\pi_{f_m}^{(J)}) -1    } \right]^2.
\end{multline}
Now we use that, for all $x >0$,
\begin{eqnarray}
\frac{1}{1\vee x}\leqslant \frac{\log(x)}{x-1}\leqslant\frac{1}{1\wedge x}.
\end{eqnarray}
Hence we infer that
$$\frac{1}{(1+\delta)^2}
\mathcal{X}_m^2\ind_{\Omega_m(\delta)}\leqslant V^2( \pi_{f_m},\pi_{\hat f_m}) \ind_{\Omega_m(\delta)}\leqslant \frac{1}{(1-\delta)^2}
\mathcal{X}_m^2\ind_{\Omega_m(\delta)},$$
with $\mathcal{X}_m^2$ defined in \eref{chideux}.
This entails that \eref{step1} is proved. It remains now to check that
$$\frac{4\rho^2 \vert m\vert}{n} \leqslant \mathbb{E}_{f_0}[ \mathcal{X}_m^2]\leqslant\frac{ 2\vert m\vert}{n}.$$
According to Lemma \ref{Projhisto} , \mbox{ for all} partition $J\in m$ and for any $ x_i \in J$,
\begin{eqnarray*}
\pi_{\hat f_m}(x_i)=\pi_{\hat f_m}^{(J)},\qquad \mbox{ with } &&\qquad
\pi_{\hat f_m}^{(J)}=\frac{1}{\vert J\vert}\sum_{i\in J}Y_i,\\
\mbox{ and }
\pi_{f_m}(x_i)=\pi_{ f_m}^{(J)}, \qquad \mbox{ with }&&\pi_{f_m}^{(J)}=\frac{1}{\vert J\vert}\sum_{i\in J}\pi_{f_0}(x_i).\end{eqnarray*}
Consequently, 
\begin{eqnarray}
\label{chideux2}
\mathcal{X}_m^2=\frac{1}{n}\sum_{J\in m} \vert J\vert \frac{(\sum_{k\in J}\varepsilon_k)^2}{\sum_{k\in J}\pi_{f_0}(x_k)[\vert J\vert-\sum_{k\in J}\pi_{f_0}(x_k)]}
=\frac{1}{n}\sum_{J\in m} \frac{(\sum_{k\in J}\varepsilon_k)^2}{\vert J\vert\pi_{f_m}^{(J)}[1-\pi_{f_m}^{(J)}]},\notag
\end{eqnarray}
and finally
\begin{eqnarray*}
\mathbb{E}_{f_0}(\mathcal{X}_m^2)=
\frac{1}n \sum_{J \in m} \mathbb{E}\left(\frac{(\sum_{k\in J}\varepsilon_k)^2}{ \vert J\vert\pi_{f_m}^{(J)} [1-\pi_{f_m}^{(J)}]     }\right)=
\frac{1 }{n}\sum_{J \in m} 
\left(\frac{1}{\vert J\vert\pi_{f_m}^{(J)}[ 1-\pi_{f_m}^{(J)}] }\right)\sum_{k\in J}\mbox{Var}\left(Y_k\right).
\end{eqnarray*}
Consequently
\begin{eqnarray*}
\mathbb{E}_{f_0}(\mathcal{X}_m^2)
=\frac{1 }{n}\sum_{J \in m} \frac{\sum_{i\in J}\pi_{f_0}(x_i)(1-\pi_{f_0}(x_i))}{\vert J\vert\pi_{f_m}^{(J)} [1-\pi_{f_m}^{(J)}]
}.
\end{eqnarray*}
Now, according to Assumption \eref{H0}, and  Lemma \ref{Projhisto}, \mbox{ for all partition }$m$, all $ J\in m$,\mbox{ and all } $x_i \in J$
\begin{eqnarray*}
 0<\rho^2 \leqslant \pi_{f_0}(x_i)(1-\pi_{f_0}(x_i)) \leqslant 1/4 , \mbox{ and } 0<\rho \leqslant \pi_{f_m}^{(J)} \mbox{ and } 0<\rho \leqslant (1-\pi_{f_m}^{(J)}). 
\end{eqnarray*}
It follows that 
\begin{eqnarray*}
4\rho^2\  \leqslant \frac{\sum_{k\in J} \pi_{f_0}(x_k) (1-\pi_{f_0}(x_k))}{ \vert J\vert\pi_{f_m}^{(J)} [1-\pi_{f_m}^{(J)}]}=\frac{\sum_{k\in J} \pi_{f_0}(x_k) (1-\pi_{f_0}(x_k))}{ \vert J\vert\pi_{f_m}^{(J)}}+\frac{\sum_{k\in J} \pi_{f_0}(x_k) (1-\pi_{f_0}(x_k))}{ \vert J\vert[1-\pi_{f_m}^{(J)}]} \leqslant 2,
\end{eqnarray*}  
and thus
\begin{eqnarray*}
\frac{4\rho^2 \vert m\vert}{n}\leqslant\frac{1 }{n}\sum_{J \in m} \frac{\sum_{i\in J}\pi_{f_0}(x_i)(1-\pi_{f_0}(x_i))}{\vert J\vert\pi_{f_m}^{(J)}[1-\pi_{f_m}^{(J)}]}
\leqslant \frac{2\vert m \vert }{n}.
\end{eqnarray*}
In other words,
\begin{eqnarray*}
\frac{4\rho^2 \vert m\vert }{n}\leqslant\mathbb{E}_{f_0}(\mathcal{X}_m^2)\leqslant \frac{2\vert m\vert}{n}.
\end{eqnarray*}
The ends up the proof of \eref{step1} and \eref{chideux}.}

\paragraph{$ \bullet$ Proof of \eref{step2} :}
We start from \eref{k_m_m},
 apply Assumption \eref{H0} and  Lemma \ref{Projhisto}, to obtain that and \eref{step2} is checked since
\begin{eqnarray*}
\vert \mathbb{E}\left(\mathcal{K}(\mathbb{P}_{f_m}^{(n)},\mathbb{P}_{\hat{f_{m}}}^{(n)})\ind_{\Omega_m^c(\delta)}\right)\vert 
&\leqslant & \frac{1}{n}\sum_{i=1}^n\mathbb{E}\left\vert \left[\log \left(
\frac{\pi_{f_m}(x_i)}{\pi_{\hat f_m}(x_i)}\right) \ind_{\Omega_m^c(\delta)}\right]\right\vert
+
\frac{1}{n}\sum_{i=1}^n\mathbb{E}
\left\vert \left[\log \left(
\frac{(1-\pi_{f_m}(x_i))}{(1-\pi_{\hat f_m}(x_i))}\right) \ind_{\Omega_m^c(\delta)}\right]\right\vert
\\
&\leqslant & 2\log\left(    \frac{1}{\rho}\right)\mathbb{P}[ \Omega_m^c(\delta)].
\end{eqnarray*}

\paragraph{$\bullet$ Proof of  \eref{step3}:}
We come to the control of $\mathbb{P}_{f_0}[ \Omega_m^c(\delta)]$.
Since 
\begin{eqnarray*}
\mathbb{P}[ \Omega_m^c(\delta)] &\leqslant& \sum_{J\in m} \mathbb{P}\left\lbrace \left\vert \frac{\pi_{\hat f_m}^{(J)}}{ \pi_{ f_m}^{(J)}  } -1\right\vert   \geqslant \delta\right\rbrace
+\sum_{J\in m} \mathbb{P}\left\lbrace\left\vert \frac{1-\pi_{\hat f_m}^{(J)}}{ 1-\pi_{ f_m}^{(J)} } -1\right\vert\geqslant 
\delta\right\rbrace,
\end{eqnarray*}
by applying  Lemma \ref{Projhisto}, we infer that
\begin{eqnarray*}
\mathbb{P}\left\lbrace \left\vert \frac{\pi_{\hat f_m}^{(J)}}{ \pi_{ f_m}^{(J)}  } -1\right\vert \geqslant 
\delta\right\rbrace &=&\mathbb{P}\left\lbrace \left\vert 
\frac{\sum_{k\in J} \varepsilon_k}{\sum_{k\in J}\pi_{ f_0}(x_k)} \right\vert \geqslant 
\delta\right\rbrace= \mathbb{P}\left\lbrace \left\vert 
\sum_{k\in J} \varepsilon_k \right\vert \geqslant 
\delta \sum_{k\in J}\pi_{ f_0}(x_k)\right\rbrace, 
\end{eqnarray*}
and 
\begin{eqnarray*}
\mathbb{P}\left\lbrace \left\vert \frac{1-\pi_{\hat f_m}^{(J)}}{ 1-\pi_{ f_m}^{(J)}  } -1\right\vert \geqslant 
\delta\right\rbrace &=&\mathbb{P}\left\lbrace \left\vert 
\frac{\sum_{k\in J} \varepsilon_k}{\sum_{k\in J}(1-\pi_{ f_0}(x_k))} \right\vert \geqslant 
\delta\right\rbrace
= \mathbb{P}\left\lbrace \left\vert 
\sum_{k\in J} \varepsilon_k \right\vert \geqslant 
\delta \sum_{k\in J}(1-\pi_{ f_0}(x_k))\right\rbrace.
\end{eqnarray*}
 We write
\begin{eqnarray*}
\mathbb{P}\left\lbrace \left\vert 
\sum_{k\in J} \varepsilon_k \right\vert \geqslant
\delta \sum_{k\in J}\pi_{ f_0}(x_k)\right\rbrace  \!\!\!&\leqslant&\!\!\!
\mathbb{P}\left\lbrace \left\vert 
\sum_{k\in J} \varepsilon_k \right\vert \geqslant
\delta \sum_{k\in J}\pi_{ f_0}(x_k) (1-\pi_{f_0}(x_k)) \right\rbrace
\end{eqnarray*}
and 
\begin{eqnarray*}
\mathbb{P}\left\lbrace \left\vert 
\sum_{k\in J} \varepsilon_k \right\vert \geqslant 
\delta \sum_{k\in J}(1-\pi_{ f_0}(x_k))\right\rbrace
\!\!\!&\leqslant&\!\!\!
\mathbb{P}\left\lbrace \left\vert 
\sum_{k\in J} \varepsilon_k \right\vert \geqslant
\delta \sum_{k\in J}\pi_{ f_0}(x_k) (1-\pi_{f_0}(x_k)) \right\rbrace.
\end{eqnarray*}
Then we have

\begin{eqnarray*}
\mathbb{P}[ \Omega_m^c(\delta)] &\leqslant& 2\sum_{J\in m}\mathbb{P}\left\lbrace \left\vert 
\sum_{k\in J} \varepsilon_k \right\vert \geqslant
\delta \sum_{k\in J}\pi_{ f_0}(x_k) (1-\pi_{f_0}(x_k)) \right\rbrace.
\end{eqnarray*}
Now, we apply Bernstein Concentration Inequality (see Massart \citeyear{massart2007} for example) to the right hand side of previous inequality, starting by recalling this Bernstein inequality.
\begin{theo}\label{Bernstein}
Let $Z_1,\cdots,Z_n$ be independent real valued random variables. Assume that there exist
some positive numbers $v$ and $c$ such that for all $k\geqslant 2$,
$$\sum_{i=1}^n\mathbb{E}\left[  \vert Z_i\vert ^k\right] \leqslant \frac{k!}{2} v c^{k-2}.$$
Then for any positive $z$,
$$\mathbb{P}\left( \sum_{i=1}^n (Z_i-\mathbb{E}(Z_i) \geqslant \sqrt{2vz}+cz\right) \leqslant \exp(-z),\mbox{ and }
\mathbb{P}\left( \sum_{i=1}^n (Z_i-\mathbb{E}(Z_i) \geqslant z\right) \leqslant \exp\left(-\frac{z^2}{2(v+cz)}\right).$$
Especially, if  $\vert Z_i \vert \leqslant b$ for all $i$, then
\begin{eqnarray}
\label{casborne}
\mathbb{P}\left( \sum_{i=1}^n (Z_i-\mathbb{E}(Z_i) \geqslant z\right) \leqslant \exp\left(-\frac{z^2}{2(\sum_{i=1}^n \mathbb{E}(Z_i^2)+bz/3)}\right).
\end{eqnarray}
\end{theo}
Applying \eref{casborne} with $z=\delta \sum_{k\in J}\pi_{ f_0}(x_k) (1-\pi_{f_0}(x_k) )$, $b=1$ and  $v=\sum_{k\in J} \pi_{ f_0}(x_k)(1-\pi_{f_0}(x_k) ),$ we get
that $$ \mathbb{P}\left\lbrace \left\vert 
\sum_{k\in J} \varepsilon_k \right\vert \geqslant
\delta \sum_{k\in J}\pi_{ f_0}(x_k) (1-\pi_{f_0}(x_k)) \right\rbrace$$
is less than
\begin{eqnarray*}
2\exp\left( -\frac{\delta^2 [  \sum_{k\in J}\pi_{ f_0}(x_k) (1-\pi_{f_0}(x_k))]^2}{2\left(  \sum_{k\in J}\pi_{ f_0}(x_k)(1-\pi_{f_0}(x_k) )  
+(\delta/3) \sum_{k\in J}\pi_{ f_0}(x_k) (1-\pi_{f_0}(x_k) )\right)}\right),
\end{eqnarray*}
and consequently
\begin{eqnarray*}
\mathbb{P}\left\lbrace \left\vert 
\sum_{k\in J} \varepsilon_k \right\vert \geqslant
\delta \sum_{k\in J}\pi_{ f_0}(x_k) (1-\pi_{f_0}(x_k)) \right\rbrace
\!\!\!&\leqslant& \!\!\! 2\exp\left[ -\frac{\delta^2}{2(1+\delta/3)} \left( \sum_{k\in J}\pi_{ f_0}(x_k) (1-\pi_{f_0}(x_k)) \right)\right]\\
\!\!\!&\leqslant&\!\!\!  2\exp\left[ -\frac{\delta^2}{2(1+\delta/3)} \vert J\vert \rho^2\right].
\end{eqnarray*}
Consequently, 
\begin{eqnarray*}
\mathbb{P}[ \Omega_m^c(\delta)] &\leqslant& 4\vert m \vert \exp(-\Delta \rho^2\Gamma [\log(n)]^2), \qquad \mbox{ with } \qquad \Delta=\frac{\delta^2}{2(1+\delta/3)},
\end{eqnarray*}
where $\Gamma$ is given by \eref{def_Gam}.
For $\epsilon>0$ and $\delta$ such that 
\begin{eqnarray}
\label{conddelta}
\frac{\delta^2}{2(1+\delta/3)}\rho^2\Gamma \log(n) \geqslant 2+\epsilon,
\end{eqnarray}
 using that  $\vert m\vert \leqslant n$ implies that 
  \begin{eqnarray*}
4\vert m\vert \exp\left(-\frac{\delta^2}{2(1+\delta/3)}\rho^2 \Gamma [\log(n)]^2 \right)
\leqslant \frac{\kappa}{n^{(1+\epsilon)}}.
\end{eqnarray*}
And Result \eref{step3} follows.

\subsection{Proof  of Theorem~\ref{theo2}}~\\

By definition, for all $m \in \mathcal{M}$,
$$\gamma_n(\hat{f}_{\hat m})+\mbox{pen}(\hat{m})\leqslant \gamma_n(\hat{f_m})+\mbox{pen}(m) \leqslant \gamma_n(f_m)+\mbox{pen}(m).$$
Applying  Formula \eref{decompgamma}, we have
\begin{equation}\label{decomp2}
 \gamma(\hat{f}_{\hat m})-\gamma(f_{0})\leqslant \gamma(f_m)-\gamma(f_{0})+\langle \vec{\varepsilon}, \hat{f}_{\hat m}-f_m\rangle_n+\mbox{pen}(m)-\mbox{pen}(\hat{m}).
\end{equation}
Following Baraud \citeyear{Baraud2000} or
Castellan \citeyear{Castellan2003},  instead of bounding the supremum of the empirical process $\langle \vec{\varepsilon}, \hat{f}_{\hat m}-f_m\rangle_n$,  we split it in three terms.
Let $$\overline{\gamma}_n(t)= {\gamma}_n(t)-\mathbb{E}_{f_0}(\gamma_n(t))=- <\vec{\varepsilon},f>_n$$
with  $<\vec{\varepsilon},f>_n $ defined in \eref{decompgamma}, and 
 write 
\begin{eqnarray*}
 \gamma(\hat{f}_{\hat m})-\gamma(f_{0})&\leqslant& \gamma(f_m)-\gamma(f_{0})+\mbox{pen}(m)-\mbox{pen}(\hat{m})\nonumber\\ &&+\overline{\gamma}_n(f_m)-\overline{\gamma}_n(f_0) 
+\overline{\gamma}_n(f_0)-\overline{\gamma}_n(f_{\hat m})+\overline{\gamma}_n(f_{\hat m})-\overline{\gamma}_n(\hat{f}_{\hat m}).
\end{eqnarray*}
In other words,
\begin{eqnarray}
 \mathcal{K}(\mathbb{P}_{f_0}^{(n)},\mathbb{P}_{\hat{f}_{\hat m}}^{(n)})&\leqslant& \mathcal{K}(\mathbb{P}_{f_0}^{(n)},\mathbb{P}_{ f_{m}}^{(n)})+\mbox{pen}(m)-\mbox{pen}(\hat{m})\nonumber \\ &&+\overline{\gamma}_n(f_m)-\overline{\gamma}_n(f_0) 
+\overline{\gamma}_n(f_0)-\overline{\gamma}_n(f_{\hat m})+\overline{\gamma}_n(f_{\hat m})-\overline{\gamma}_n(\hat{f}_{\hat m}).\label{decomp3}
\end{eqnarray}
The proof of Theorem~\ref{theo2} can be decomposed in three steps :
\begin{enumerate}
\item \label{biais1} We prove that for $\epsilon > 0,$ 
$$\mathbb{E}_{f_0}\big[(\overline{\gamma}_n(f_m)-\overline{\gamma}_n(f_0))\ind_{\Omega_{m_{f}}(\delta)}\big]
\leqslant \frac{\kappa^{\prime}(\rho,\delta,\Gamma,\epsilon)}{n^{(1+\epsilon)}}.$$

\item \label{deviation} 
Let $\Omega_1(\xi)$ be the event
\begin{eqnarray*}
\Omega_1(\xi)=\bigcap_{m^\prime\in \mathcal{M}}\left\{ \chi_{m^\prime}^{2}\ind_{\Omega_{m_{f}}(\delta)}\right.
\!\!\!&\leqslant& \!\!\!\left.
\frac{2}{n}\vert m^\prime\vert+\frac{16}{n}\Big(1+\frac{\delta}{3}\Big)\sqrt{(L_{m^{\prime}}\vert m^{\prime}\vert+\xi)\vert m^\prime\vert}+\frac{8}{n}\Big(1+\frac{\delta}{3}\Big)(L_{m^{\prime}}\vert m^{\prime}\vert+\xi)\right\},\end{eqnarray*}
where $(L_{m^{\prime}})_{m'\in{\mathcal M}}$ satisfies Condition (\ref{sigma}) and $m_{f}$ is given by Definition~\ref{def_mf}. For all $m^{\prime}$ in $\mathcal{M}$ we prove that on  $\Omega_1(\xi)$
   \begin{eqnarray}\label{chy2}
  \notag\Big( \overline{\gamma}_n(f_{m^{\prime}})-\overline{\gamma}_n(\hat{f}_{m^{\prime}})\Big)\ind_{\Omega_{m_{f}}(\delta)} \leqslant
  &&\frac{1}{2n}\Big(\frac{1+\delta}{1-\delta}\Big)\vert m^{\prime}\vert \Big[2+\Big(1+\frac{\delta}{3}\Big)\Big(2\delta +8L_{m^\prime}+16\sqrt{L_{m^\prime}}\Big) \Big]\\
&&+\frac{4\xi}{n} \Big(\frac{1+\delta}{1-\delta}\Big)\Big(1+\frac{\delta}{3}\Big)\Big(1+\frac{4}{\delta}\Big)+\frac{1}{1+\delta} \mathcal{K}(\mathbb{P}_{f_{m^{\prime}}}^{(n)},\mathbb{P}_{\hat{f}_{m^{\prime}}}^{(n)})\ind_{\Omega_{m_{f}}(\delta)},
 \end{eqnarray}
 and
 \begin{equation}\label{probchy2}
  \mathbb{P}(\Omega_1(\xi)^{c}) \leqslant 2\Sigma e^{-\xi}.
\end{equation}

\item  \label{biais2} 
Let $\Omega_2(\xi)$ be the event
$$\Omega_2(\xi)=\bigcap_{m^\prime\in \mathcal{M}}\left[(\overline{\gamma}_n(f_0)-\overline{\gamma}_n(f_{m^{\prime}}))  \leqslant
\mathcal{K}(\mathbb{P}_{f_0}^{(n)},\mathbb{P}_{f_{m^\prime}}^{(n)})-2 h^2(\mathbb{P}_{f_0}^{(n)}, \mathbb{P}_{f_{m^\prime}}^{(n)})+\frac{2}{n}(L_m^{\prime}\vert m^{\prime}\vert+\xi)\right].$$
We prove that,
$\mathbb{P}(\Omega_2(\xi)^{c}) \leqslant \Sigma e^{-\xi}.$
\end{enumerate}
Now, we will prove the result of Theorem~\ref{theo2} using (R-\ref{biais1}), (R-\ref{deviation}) and (R-\ref{biais2}).\\
According to \eref{decomp3},  we can write 
\begin{eqnarray*}
 \mathcal{K}(\mathbb{P}_{f_0}^{(n)},\mathbb{P}_{\hat{f}_{\hat m}}^{(n)}) \ind_{\Omega_{m_f}(\delta)}&\leqslant& \mathcal{K}(\mathbb{P}_{f_0}^{(n)},\mathbb{P}_{f_m}^{(n)})+\mbox{pen}(m)-\mbox{pen}(\hat{m})\\
&&+(\overline{\gamma}_n(f_m)-\overline{\gamma}_n(f_0)) \ind_{\Omega_{m_f}(\delta)}+(\overline{\gamma}_n(f_0)-\overline{\gamma}_n(f_{\hat m}))\ind_{\Omega_{m_f}(\delta)}+(\overline{\gamma}_n(f_{\hat m})-\overline{\gamma}_n(\hat{f}_{\hat m}) \ind_{\Omega_{m_f}(\delta)}.
\end{eqnarray*}
Combining (R-\ref{deviation}) and (R-\ref{biais2}) with $m^{\prime}=\hat{m}$, we infer that on $\Omega_1(\xi)\bigcap \Omega_2(\xi)$  
\begin{eqnarray*}
 \mathcal{K}(\mathbb{P}_{f_0}^{(n)},\mathbb{P}_{\hat{f}_{\hat m}}^{(n)}) \ind_{\Omega_{m_f}(\delta)}&\leqslant& \mathcal{K}(\mathbb{P}_{f_0}^{(n)},\mathbb{P}_{f_m}^{(n)})+\mbox{pen}(m)-\mbox{pen}(\hat{m})+(\overline{\gamma}_n(f_m)-\overline{\gamma}_n(f_0)) \ind_{\Omega_{m_f}(\delta)}\\
 &&+\frac{1}{2n}\Big(\frac{1+\delta}{1-\delta}\Big)\vert\hat{m}\vert \Big[2+\Big(1+\frac{\delta}{3}\Big)\Big(2\delta +8L_{\hat{m}}+16\sqrt{L_{\hat{m}}}\Big) \Big]+2L_{\hat{m}}\frac{\vert \hat{m}\vert}{n}\\
&&+\frac{4\xi}{n}\Big[\frac{1}{2}+\Big(\frac{1+\delta}{1-\delta}\Big)\Big(1+\frac{\delta}{3}\Big)\Big(1+\frac{4}{\delta}\Big)\Big]\\
&&+\Big[\mathcal{K}(\mathbb{P}_{f_0}^{(n)},\mathbb{P}_{f_{\hat{m}}}^{(n)})-2 h^2(\mathbb{P}_{f_0}^{(n)}, \mathbb{P}_{f_{\hat{m}}}^{(n)})+\frac{1}{1+\delta} \mathcal{K}(\mathbb{P}_{f_{\hat{m}}}^{(n)},\mathbb{P}_{\hat{f}_{\hat{m}}}^{(n)})\Big]\ind_{\Omega_{m_{f}}(\delta)}.
\end{eqnarray*}
This implies that
\begin{eqnarray*}
 \mathcal{K}(\mathbb{P}_{f_0}^{(n)},\mathbb{P}_{\hat{f}_{\hat m}}^{(n)}) \ind_{\Omega_{m_f}(\delta)}&\leqslant& \mathcal{K}(\mathbb{P}_{f_0}^{(n)},\mathbb{P}_{f_m}^{(n)})+\mbox{pen}(m)-\mbox{pen}(\hat{m})+(\overline{\gamma}_n(f_m)-\overline{\gamma}_n(f_0)) \ind_{\Omega_{m_f}(\delta)}\\
 && + \frac{\vert\hat{m}\vert}{n}\Big[\Big(\frac{1+\delta}{1-\delta}\Big)+\Big(\frac{\delta(1+\delta)^{2}}{1-\delta}\Big)
 +\Big(\frac{(1+\delta)^{2}}{1-\delta}\Big)\Big(6L_{\hat{m}}+8\sqrt{L_{\hat{m}}}\Big)\Big]\\
&&+\frac{4\xi}{n}\Big[\frac{1}{2}+\Big(\frac{1+\delta}{1-\delta}\Big)\Big(1+\frac{\delta}{3}\Big)\Big(1+\frac{4}{\delta}\Big)\Big]\\
&&+\Big[\mathcal{K}(\mathbb{P}_{f_0}^{(n)},\mathbb{P}_{f_{\hat{m}}}^{(n)})-2 h^2(\mathbb{P}_{f_0}^{(n)}, \mathbb{P}_{f_{\hat{m}}}^{(n)}))+\frac{1}{1+\delta} \mathcal{K}(\mathbb{P}_{f_{\hat{m}}}^{(n)},\mathbb{P}_{\hat{f}_{\hat{m}}}^{(n)})\Big]\ind_{\Omega_{m_{f}}(\delta)}.
\end{eqnarray*}
Since $$\Big\{\Big(\frac{1+\delta}{1-\delta}\Big)(1+\delta(1+\delta))\vee \Big(\frac{(1+\delta)^{2}}{1-\delta}\Big)\Big\} \leqslant  C(\delta) \mbox{ with }
C(\delta):=\Big(\frac{1+\delta}{1-\delta}\Big)^{3},$$
we infer
\begin{eqnarray*}
 \mathcal{K}(\mathbb{P}_{f_0}^{(n)},\mathbb{P}_{\hat{f}_{\hat m}}^{(n)}) \ind_{\Omega_{m_f}(\delta)}&\leqslant& \mathcal{K}(\mathbb{P}_{f}^{(n)},\mathbb{P}_{f_m}^{(n)})+\mbox{pen}(m)-\mbox{pen}(\hat{m})+(\overline{\gamma}_n(f_m)-\overline{\gamma}_n(f_0)) \ind_{\Omega_{m_f}(\delta)}\\
 && +\frac{\vert\hat{m}\vert}{n}C(\delta) \Big[1+6L_{\hat{m}}+8\sqrt{L_{\hat{m}}}\Big]+\frac{4\xi}{n}\Big[\frac{1}{2}+\Big(\frac{1+\delta}{1-\delta}\Big)\Big(1+\frac{\delta}{3}\Big)\Big(1+\frac{4}{\delta}\Big)\Big]\\
&&+\Big[\mathcal{K}(\mathbb{P}_{f_0}^{(n)},\mathbb{P}_{f_{\hat{m}}}^{(n)})-2 h^2(\mathbb{P}_{f_0}^{(n)}, \mathbb{P}_{f_{\hat{m}}}^{(n)})+\frac{1}{1+\delta} \mathcal{K}(\mathbb{P}_{f_{\hat{m}}}^{(n)},\mathbb{P}_{\hat{f}_{\hat{m}}}^{(n)})\Big]\ind_{\Omega_{m_{f}}(\delta)}.
\end{eqnarray*}
Using  Pythagore's type identity $\mathcal{K}(\mathbb{P}_{f_0},\mathbb{P}_{\hat{f}_{\hat m}})=\mathcal{K}(\mathbb{P}_{f_0}^{(n)},\mathbb{P}_{f_{\hat{m}}}^{(n)})+\mathcal{K}(\mathbb{P}_{f_{\hat{m}}}^{(n)},\mathbb{P}_{\hat{f}_{\hat m}}^{(n)})$ (see Equation (7.42) in Massart \citeyear{massart2007}) we have

\begin{eqnarray*}
 \mathcal{K}(\mathbb{P}_{f_0}^{(n)},\mathbb{P}_{\hat{f}_{\hat m}}^{(n)}) \ind_{\Omega_{m_f}(\delta)}&\leqslant& \mathcal{K}(\mathbb{P}_{f}^{(n)},\mathbb{P}_{f_m}^{(n)})+\mbox{pen}(m)-\mbox{pen}(\hat{m})+(\overline{\gamma}_n(f_m)-\overline{\gamma}_n(f_0)) \ind_{\Omega_{m_f}(\delta)}\\
 && +\frac{\vert\hat{m}\vert}{n}C(\delta) \Big[1+6L_{\hat{m}}+8\sqrt{L_{\hat{m}}}\Big]+\frac{4\xi}{n}\Big[\frac{1}{2}+\Big(\frac{1+\delta}{1-\delta}\Big)\Big(1+\frac{\delta}{3}\Big)\Big(1+\frac{4}{\delta}\Big)\Big]\\
&&+\Big[\mathcal{K}(\mathbb{P}_{f_0}^{(n)},\mathbb{P}_{\hat f_{\hat{m}}}^{(n)})-2 h^2(\mathbb{P}_{f_0}^{(n)}, \mathbb{P}_{f_{\hat{m}}}^{(n)})-\frac{\delta}{1+\delta} \mathcal{K}(\mathbb{P}_{f_{\hat{m}}}^{(n)},\mathbb{P}_{\hat{f}_{\hat{m}}}^{(n)})\Big]\ind_{\Omega_{m_{f}}(\delta)}.
\end{eqnarray*}
Now, we successively use
\begin{itemize}
\item[(i)] the relation between Kullback-Leibler information and the Hellinger distance $ \mathcal{K}(\mathbb{P}_{f_{\hat{m}}}^{(n)},\mathbb{P}_{\hat{f}_{\hat{m}}}^{(n)})\geq 2 h^2(\mathbb{P}_{f_{\hat{m}}}^{(n)},\mathbb{P}_{\hat{f}_{\hat{m}}}^{(n)})$ (see Lemma 7.23 in Massart \citeyear{massart2007}),
\item [(ii)] and inequality $ h^2(\mathbb{P}_{f_0}^{(n)}, \mathbb{P}_{\hat{f}_{\hat{m}}}^{(n)})\leqslant2[h^2(\mathbb{P}_{f_0}^{(n)}, \mathbb{P}_{f_{\hat{m}}}^{(n)})+  h^2(\mathbb{P}_{f_{\hat{m}}}^{(n)}, \mathbb{P}_{\hat{f}_{\hat{m}}}^{(n)})]$.
\end{itemize}
Consequently, on $\Omega_1(\xi)\bigcap \Omega_2(\xi)$ 
\begin{eqnarray*}
 \frac{\delta}{1+\delta}h^{2}(\mathbb{P}_{f_0}^{(n)},\mathbb{P}_{\hat{f}_{\hat m}}^{(n)}) \ind_{\Omega_{m_f}(\delta)}&\leqslant& \mathcal{K}(\mathbb{P}_{f_0}^{(n)},\mathbb{P}_{f_m}^{(n)})+\mbox{pen}(m)-\mbox{pen}(\hat{m})+(\overline{\gamma}_n(f_m)-\overline{\gamma}_n(f_0)) \ind_{\Omega_{m_f}(\delta)}\\
 && +\frac{\vert\hat{m}\vert }{n} C(\delta) \Big[1+6L_{\hat{m}}+8\sqrt{L_{\hat{m}}}\Big]+\frac{4\xi}{n} \Big[\frac{1}{2}+\Big(\frac{1+\delta}{1-\delta}\Big)\Big(1+\frac{\delta}{3}\Big)\Big(1+\frac{4}{\delta}\Big)\Big].
\end{eqnarray*}
Since $\mbox{pen}(\hat{m})\geq \mu \vert \hat{m}\vert\Big[1+6L_{\hat{m}}+8\sqrt{L_{\hat{m}}}\Big]/n$, by taking $\mu=C(\delta)$ yields
that on $\Omega_1(\xi)\bigcap \Omega_2(\xi)$
\begin{eqnarray*}
 h^{2}(\mathbb{P}_{f_0},\mathbb{P}_{\hat{f}_{\hat m}}) \ind_{\Omega_{m_f}(\delta)}&\leqslant& 
\frac{2\mu^{1/3}}{\mu^{1/3}-1}\Big( \mathcal{K}(\mathbb{P}_{f_0}^{(n)},\mathbb{P}_{f_m}^{(n)})+\mbox{pen}(m)+(\overline{\gamma}_n(f_m)-\overline{\gamma}_n(f_0)) \ind_{\Omega_{m_f}(\delta)}\Big)+ \frac{\xi}{n} C_1(\mu).
\end{eqnarray*}
Then, using that
$$\mathbb{P}(\Omega_1(\xi)^{c}\cup \Omega_2(\xi)^{c}) \leqslant 3\Sigma e^{-\xi},$$
we deduce that $\mathbb{P}(\Omega_1(\xi)\cap \Omega_2(\xi))\geq 1-3\Sigma e^{-\xi}.$
We now integrating with respect to $\xi$, and use (R-\ref{biais1}) to write that
\begin{eqnarray*}
 \mathbb{E}_{f_0}\Big[h^{2}(\mathbb{P}_{f_0},\mathbb{P}_{\hat{f}_{\hat m}}) \ind_{\Omega_{m_f}(\delta)}\Big]&\leqslant& \frac{2\mu^{1/3}}{\mu^{1/3}-1}\Big( \mathcal{K}(\mathbb{P}_{f_0}^{(n)},\mathbb{P}_{f_m}^{(n)})+\mbox{pen}(m)\Big)+ \frac{\kappa_1(\rho,\mu,\Gamma,\epsilon)}{n^{(1+\epsilon)}}+\frac{C_2(\mu,\Sigma)}{n} .
\end{eqnarray*}
Furthermore, since $h^{2}(\mathbb{P}_{f_0},\mathbb{P}_{\hat{f}_{\hat m}})\leqslant1,$ by applying Inequality (\ref{step3}) we have,
$$\mathbb{E}_{f_0}\Big[h^{2}(\mathbb{P}_{f_0},\mathbb{P}_{\hat{f}_{\hat m}}) \ind_{\Omega_{m_f}^{c}(\delta)}\Big]\leq\frac{\kappa_2(\rho,\mu,\Gamma,\epsilon)}{n^{(1+\epsilon)}}. $$
Hence we conclude that 
\begin{eqnarray*}
 \mathbb{E}_{f_0}\Big[h^{2}(\mathbb{P}_{f_0},\mathbb{P}_{\hat{f}_{\hat m}})\Big]&\leqslant& \frac{2\mu^{1/3}}{\mu^{1/3}-1}\Big( \mathcal{K}(\mathbb{P}_{f_0}^{(n)},\mathbb{P}_{f_m}^{(n)})+\mbox{pen}(m)\Big)+ \frac{\kappa_3(\rho,\mu,\Gamma,\epsilon)}{n^{(1+\epsilon)}}+\frac{C_2(\mu,\Sigma)}{n} ,
\end{eqnarray*}
 and minimizing over $\mathcal{M}$ leads to the result of Theorem~\ref{theo2}.\\
 We now come to the proofs of (R-\ref{biais1}), (R-\ref{deviation}) and (R-\ref{biais2}).\\
 $\bullet$ Proof of  (R-\ref{biais1})\\ 
We know that
\begin{eqnarray*}
\Big\vert  \mathbb{E}_{f_0}\Big[(\overline{\gamma}_n(f_m)-\overline{\gamma}_n(f_0))\ind_{\Omega_{m_f}(\delta)}\Big] \Big\vert\!\!\!&=&\!\!\!
\Big\vert  \mathbb{E}_{f_0}\Big[(\overline{\gamma}_n(f_m)-\overline{\gamma}_n(f_0))\ind_{\Omega_{m_f}^{c}(\delta)}\Big] \Big\vert\\
\!\!\!&\leq&\!\!\!  \mathbb{E}_{f_0}\Big[\frac{1}{n}\sum_{i=1}^{n}\Big\{\Big\vert\epsilon_{i}\log\{\frac{\pi_{f_{m}}(x_i)}{\pi_{f_{0}}(x_i)}\} \Big\vert + \Big\vert \epsilon_i \log\{\frac{1-\pi_{f_m}(x_i)}{1-\pi_{f_0}(x_i)}\}\Big\vert\Big\}\ind_{\Omega_{m_f}^{c}(\delta)}\Big]\\
\!\!\!&\leq&\!\!\! 2\log\left\lbrace\frac{1}{\rho}\right\rbrace\mathbb{P}(\Omega_{m_f}^{c}(\delta)).
\end{eqnarray*}
We conclude the proof of (R-\ref{biais1}) by
using Inequality (\ref{step3}), which  implies that
$$\Big\vert  \mathbb{E}_{f_0}\Big[(\overline{\gamma}_n(f_m)-\overline{\gamma}_n(f_0))\ind_{\Omega_{m_f}(\delta)}\Big] \Big\vert\leq
2\log\left\lbrace\frac{1}{\rho}\right\rbrace\frac{\kappa(\rho,\delta,\Gamma,\epsilon)}{n^{(1+\epsilon)}}=\frac{\kappa^{\prime}(\rho,\delta,\Gamma,\epsilon)}{n^{(1+\epsilon)}}.$$

$\bullet$ Proof of  (R-\ref{deviation})\\
We start by the proof of (\ref{chy2})
\begin{eqnarray*}
\overline{\gamma}_n(f_{m^{\prime}})-\overline{\gamma}_n(\hat{f}_{m^{\prime}})
&=& -\frac{1}{n}\sum_{i=1}^{n}\Big\{\epsilon_{i}\log\Big(\frac{\pi_{f_{m^{\prime}}}(x_i)}{\pi_{\hat{f}_{m^{\prime}}}(x_i)}\Big)  - \epsilon_i \log\Big(\frac{1-\pi_{f_{m^{\prime}}}(x_i)}{1-\pi_{\hat{f}_{m^{\prime}}}(x_i)}\Big)\Big\}\\
&=& -\frac{1}{n}\sum_{J\in m^{\prime}}\Big(\sum_{i\in J}\epsilon_{i}\Big)\Big[\frac{\sqrt{\vert J \vert\pi_{f_{m^{\prime}}}^{(J)}}}{\sqrt{\vert J \vert\pi_{f_{m^{\prime}}}^{(J)}}}\log\Big(\frac{\pi_{f_{m^{\prime}}}^{(J)}}{\pi_{\hat{f}_{m^{\prime}}}^{(J)}}\Big)  -  
\frac{\sqrt{\vert J \vert 1-\pi_{f_{m^{\prime}}}^{(J)}}}{\sqrt{\vert J \vert (1-\pi_{f_{m^{\prime}}}^{(J)})}}
\log\Big(\frac{1-\pi_{f_{m^{\prime}}}^{(J)}}{1-\pi_{\hat{f}_{m^{\prime}}}^{(J)}}\Big)\Big].
\end{eqnarray*}
By Cauchy-Schwarz inequality, we have
\begin{multline*}
\overline{\gamma}_n(f_{m^{\prime}})-\overline{\gamma}_n(\hat{f}_{m^{\prime}})
\leq 
\sqrt{\frac{1}{n}\sum_{J\in m^{\prime}}\vert J\vert \Big[ \pi_{f_{m^{\prime}}}^{(J)}\log^{2}{\Big(\frac{\pi_{\hat{f}_{m^{\prime}}}^{(J)}}{\pi_{f_{m^{\prime}}}^{(J)}}\Big)}+(1-\pi_{f_{m^{\prime}}}^{(J)})\log^{2}\Big({\frac{1-\pi_{\hat{f}_{m^{\prime}}}^{(J)}}{1-\pi_{f_{m^{\prime}}}^{(J)}}\Big)}\Big]}\\
\\
\times\sqrt{\frac{1}{n}\sum_{J \in m^{\prime}}\Big[\frac{\Big(\sum_{i\in J}\epsilon_{i}\Big)^{2}}{\vert J \vert\pi_{f_{m^{\prime}}}^{(J)}}
+\frac{\Big(\sum_{i\in J}\epsilon_{i}\Big)^{2}}{\vert J \vert(1-\pi_{f_{m^{\prime}}}^{(J)})}\Big]} ~ ,
\end{multline*}
and in other words
\begin{eqnarray*}\overline{\gamma}_n(f_{m^{\prime}})-\overline{\gamma}_n(\hat{f}_{m^{\prime}})
\leq \sqrt{\mathcal{X}^{2}_{m^{\prime}}}\times \sqrt{V^{2}(\pi_{f_{m^{\prime}}},\pi_{\hat{f}_{m^{\prime}}})},
\end{eqnarray*}
where   $\mathcal{X}^{2}_{m^{\prime}}$  and $V^{2}(\pi_{f_{m^{\prime}}},\pi_{\hat{f}_{m^{\prime}}})$  are defined respectively in ~\eref{chideux} and \eref{V}  .
Using both that inequality  $2xy\leqslant\theta x^{2}+ \theta^{-1}y^{2}$, for all $ x>0$, $y>0$ with $\theta=(1+\delta)/(1-\delta),$ and Inequality~(\ref{encadrK}), we obtain on $\Omega_{m_{f}}(\delta)$  that,
\begin{eqnarray*}
\overline{\gamma}_n(f_{m^{\prime}})-\overline{\gamma}_n(\hat{f}_{m^{\prime}}))
\!\!\!&\leq&\!\!\! \frac{1}{2}\Big(\frac{1+\delta}{1-\delta}\Big) \chi^{2}_{m^{\prime}}+\frac{1}{1+\delta} \mathcal{K}(\mathbb{P}_{f_{m^{\prime}}}^{(n)},\mathbb{P}_{\hat{f}_{m^{\prime}}}^{(n)}). 
\end{eqnarray*}
Consequently, on $\Omega_{1}(\xi)$ 
\begin{eqnarray*}
(\overline{\gamma}_n(f_{m^{\prime}})-\overline{\gamma}_n(\hat{f}_{m^{\prime}}))\ind_{\Omega_{m_{f}}(\delta)}
&\leq&\frac{1}{2n}\Big(\frac{1+\delta}{1-\delta}\Big) \Big[2\vert m^\prime\vert+16\Big(1+\frac{\delta}{3}\Big)\sqrt{(L_{m^{\prime}}\vert m^{\prime}\vert+\xi)\vert m^\prime\vert}+8\Big(1+\frac{\delta}{3}\Big)(L_{m^{\prime}}\vert m^{\prime}\vert+\xi)\Big]\\
&+&\frac{1}{1+\delta} \mathcal{K}(\mathbb{P}_{f_{m^{\prime}}}^{(n)},\mathbb{P}_{\hat{f}_{m^{\prime}}}^{(n)})\ind_{\Omega_{m_{f}}(\delta)}.
\end{eqnarray*}
 Using inequalities $\vert x+y\vert^{1/2}\leqslant\vert x\vert^{1/2}+ \vert y \vert^{1/2}$ and $2xy\leqslant\theta x^{2}+ \theta^{-1}y^{2}$ with 
$\theta=\delta/4$, we infer that \eref{chy2} follows since
\begin{eqnarray*}
\overline{\gamma}_n(f_{m^{\prime}})-\overline{\gamma}_n(\hat{f}_{m^{\prime}}))\ind_{\Omega_{m_{f}}(\delta)}
\!\!\!&\leq&\!\!\!\frac{1}{2n}\Big(\frac{1+\delta}{1-\delta}\Big) \Big[2\vert m^\prime\vert+\Big(1+\frac{\delta}{3}\Big)\Big(16\sqrt{L_{m^{\prime}}}\vert m^{\prime}\vert+8L_{m^{\prime}}\vert m^\prime\vert+2\delta\vert m^{\prime}\vert\Big)\\ \!\!\!&&\!\!\! +8\xi\Big(1+\frac{\delta}{3}\Big)(1+\frac{4}{\delta})\Big]+\frac{1}{1+\delta} \mathcal{K}(\mathbb{P}_{f_{m^{\prime}}}^{(n)},\mathbb{P}_{\hat{f}_{m^{\prime}}}^{(n)})\ind_{\Omega_{m_{f}}(\delta)}\\
 \!\!\!&\leq&\frac{1}{2n}\Big(\frac{1+\delta}{1-\delta}\Big)\vert m^{\prime}\vert \Big[2+\Big(1+\frac{\delta}{3}\Big)\Big(2\delta +8L_{m^\prime}+16\sqrt{L_{m^\prime}}\Big) \Big]\\
\!\!\!&&+\frac{ 4\xi}{n}\Big(\frac{1+\delta}{1-\delta}\Big)\Big(1+\frac{\delta}{3}\Big)\Big(1+\frac{4}{\delta}\Big)+\frac{1}{1+\delta} \mathcal{K}(\mathbb{P}_{f_{m^{\prime}}}^{(n)},\mathbb{P}_{\hat{f}_{m^{\prime}}}^{(n)})\ind_{\Omega_{m_{f}}(\delta)}.
\end{eqnarray*}

$\bullet$ Proof of \eref{probchy2} : \\
Write 
$\mathcal{X}_{m^{\prime}}^2=\sum_{J\in m^{\prime}}\{Z_{1,J}+Z_{2,J}\},$
where $$Z_{1,J}=\frac{1}{n}\frac{(\sum_{k\in J}\varepsilon_k)^2}{\vert J\vert\pi_{f_{m^{\prime}}}^{(J)}}~~ \mbox{and}~~ Z_{2,J}=\frac{1}{n}\frac{(\sum_{k\in J}\varepsilon_k)^2}{\vert J\vert (1-\pi_{f_{m^{\prime}}}^{(J)})}.$$
We will control $\sum_{J\in m^{\prime}}Z_{1,J}$ and $\sum_{J\in m^{\prime}}Z_{2,J}$ separately.
In order to use Bernstein inequality (see Theorem~\ref{Bernstein}), we need an upper bound of  $\sum_{J\in m^{\prime}}\mathbb{E}[Z_{1,J}^{p}\ind_{\Omega_{m_f}(\delta)}]$, for every $p\geq 2$. By definition
\begin{eqnarray*}
\mathbb{E}[Z_{1,J}^{p}\ind_{\Omega_{m_f}(\delta)}]&=&\frac{1}{\Big(n\vert J\vert\pi_{f_{m^{\prime}}}^{(J)}\Big)^{p}}\int_{0}^{\infty}2px^{2p-1}
\mathbb{P}\Big(\Big\{\vert\sum_{k\in J}\varepsilon_{k}\vert\geq x\Big\}\cap \Omega_{m_f}(\delta)\Big)dx.
\end{eqnarray*}
For every $m^{\prime}$ constructed on the grid $m_f$, 
for all $ J\in m^{\prime}$,  on $\Omega_{m_f}(\delta)\cap \Big\{x\leqslant\vert\sum_{k\in J}\varepsilon_{k}\vert\Big\} ,$  we have $$x\leqslant\vert\sum_{k\in J}\varepsilon_{k}\vert\leqslant\delta \sum_{i\in J}\pi_{f_0}(x_i).$$
Combining the previous inequality, the Bernstein inequality \eref{casborne} with the fact that  $\varepsilon_k\leqslant1 $, we infer that
\begin{eqnarray*}
\mathbb{E}[Z_{1,J}^{p}\ind_{\Omega_{m_f}(\delta)}]&\leq&\frac{1}{\Big(n \sum_{k\in J}\pi_{f_{0}}(x_k)\Big)^{p}}\int_{0}^{\delta \sum_{k\in J}\pi_{f_0}(x_k)}2px^{2p-1}
\mathbb{P}\Big(\vert\sum_{k\in J}\varepsilon_{k}\vert\geq x\Big)dx\\
&\leq& \frac{1}{\Big(n \sum_{k\in J}\pi_{f_{0}}(x_k)\Big)^{p}}\int_{0}^{\delta \sum_{i\in J}\pi_{f_0}(x_i)}
4px^{2p-1}\exp\Big(-\frac{x^{2}}{2\Big( \frac{x}{3}+\sum_{k\in J}\pi_{f_{0}}(x_k) \Big)}\Big)dx\\
&\leq&  \frac{1}{\Big(n \sum_{k\in J}\pi_{f_{0}}(x_k)\Big)^{p}}\int_{0}^{\delta \sum_{i\in J}\pi_{f_0}(x_i)}
4px^{2p-1}\exp\Big(-\frac{x^{2}}{2\Big( 1+\frac{\delta}{3}\Big)\sum_{k\in J}\pi_{f_{0}}(x_k)}\Big)dx\\
&\leq&\frac{1}{n^{p}}2^{p+1}(1+\frac{\delta}{3})^{p}p\int_{0}^{\infty}t^{p-1}\exp(-t)dt\\
&\leq& \frac{1}{n^{p}}2^{p+1}p(1+\frac{\delta}{3})^{p}(p!).
\end{eqnarray*}
Consequently
$$ \sum_{J\in m^{\prime}}\mathbb{E}[Z_{1,J}^{p}\ind_{\Omega_{m_f}(\delta)}]\leqslant\frac{1}{n^{p}}2^{p+1}p(1+\frac{\delta}{3})^{p}(p!)\times\vert m^{\prime}\vert .$$
Now, since $p\leqslant2^{p-1}$, we have
$$ \sum_{J\in m^{\prime}}\mathbb{E}[Z_{1,J}^{p}\ind_{\Omega_{m_f}(\delta)}]\leqslant
\frac{p!}{2}\times\Big[\frac{32}{n^{2}}(1+\frac{\delta}{3})^{2}\vert m^{\prime}\vert \Big]\times\Big[\frac{4}{n}(1+\frac{\delta}{3})\Big]^{p-2} .$$
Using Bernstein inequality and that $\mathbb{E}\Big[\sum_{J\in m^{\prime}}Z_{1,J})\Big]\leqslant\vert m^{\prime}\vert/n$, we have that for every positive $x$
$$\mathbb{P}\Big(\sum_{J\in m^{\prime}}Z_{1,J}\ind_{\Omega_{m_f}(\delta)}\geq \frac{\vert m^{\prime}\vert}{n} + \frac{8}{n}(1+\frac{\delta}{3})\sqrt{x\vert m^{\prime}\vert}+ \frac{4}{n}(1+\frac{\delta}{3})x\Big)\leqslant\exp(-x).$$	
In the same way we prove that
$$\mathbb{P}\Big(\sum_{J\in m^{\prime}}Z_{2,J}\ind_{\Omega_{m_f}(\delta)}\geq \frac{\vert m^{\prime}\vert}{n} + \frac{8}{n}(1+\frac{\delta}{3})\sqrt{x\vert m^{\prime}\vert}+ \frac{4}{n}(1+\frac{\delta}{3})x\Big)\leqslant\exp(-x).$$
Hence
$$\mathbb{P}\Big(\mathcal{X}_{m^{\prime}}^{2}\ind_{\Omega_{m_f}(\delta)}\geq \frac{2\vert m^{\prime}\vert}{n} + \frac{16}{n}(1+\frac{\delta}{3})\sqrt{x\vert m^{\prime}\vert}+ \frac{8}{n}(1+\frac{\delta}{3})x\Big)\leqslant2\exp(-x),$$
and we conclude that
$ \mathbb{P}(\Omega_1^{c}(\xi))\leqslant2\sum_{m^{\prime}}\exp(-L_m^{\prime}\vert m^{\prime}\vert-\xi)=2\Sigma e^{-\xi}.$
This ends the proof of (R-\ref{deviation}).

$\bullet$ Proof of (R-\ref{biais2})\\
Recall that $\overline{\gamma}_n(f)=\gamma_n(f)-\mathbb{E}(\gamma_n(f))$ for every $f$. According to  Markov inequality, for $b>0$,
\begin{eqnarray*}
\mathbb{P}((\overline{\gamma}_n(f_{0})-\overline{\gamma}_n(g))\geq b)
&=&\mathbb{P}\Big(\exp\Big(\frac{n}{2}(\overline{\gamma}_n(f_{0})-\overline{\gamma}_n(g))\Big)\geq \exp\Big(\frac{nb}{2}\Big)\Big)\\
&\leq& \exp\Big(\frac{-nb}{2}\Big)\mathbb{E}\Big[\exp\Big(\frac{n}{2}(\overline{\gamma}_n(f_{0})-\overline{\gamma}_n(g))\Big)\Big]\\
&=& \exp\Big[\frac{-nb}{2}+\log\mathbb{E}\Big[\exp\Big(\frac{n}{2}\Big(\gamma_n(f_{0})-\gamma_n(g)\Big)+\frac{n}{2}\mathbb{E}\Big[\gamma_n(g)-\gamma_n(f_{0}) \Big]\Big)\Big]\\
&\leq& \exp\Big[\frac{-nb}{2}+\frac{n}{2}\mathcal{K}(\mathbb{P}_{f_0}^{(n)},\mathbb{P}_{g}^{(n)})+\log\mathbb{E}\Big[\exp\Big(\frac{n}{2}\Big(\gamma_n(f_{0})-\gamma_n(g)\Big)\Big)\Big]\Big].
\end{eqnarray*}
Now,
\begin{eqnarray*}
\log\mathbb{E}\Big[\exp\Big(\frac{n}{2}\Big(\gamma_n(f_{0})-\gamma_n(g)\Big)\Big)\Big]&=& \log \mathbb{E}\Big[\exp\Big( \frac{1}{2}\sum_{i=1}^{n}Y_i\log(\frac{\pi_{g}(x_i)}{\pi_{f_0}(x_i)})+(1-Y_i)\log(\frac{1-\pi_{g}(x_i)}{1-\pi_{f_0}(x_i)})\Big) \Big]\\
&=& \log \mathbb{E}\Big[\Pi_{i=1}^{n}\Big\{\Big(\frac{\pi_{g}(x_i)}{\pi_{f_0}(x_i)}\Big)^{Y_i/2}\times\Big(\frac{1-\pi_{g}(x_i)}{1-\pi_{f_0}(x_i)}\Big)^{(1-Y_i)/2}\Big\} \Big]\\
&=& \log\Pi_{i=1}^{n}\Big\{\sqrt{\frac{\pi_{g}(x_i)}{\pi_{f_0}(x_i)}}\pi_{f_0}(x_i)+\sqrt{\frac{1-\pi_{g}(x_i)}{1-\pi_{f_0}(x_i)}}(1-\pi_{f_0}(x_i))\Big\} \\
&=&\sum_{i=1}^{n}\log\Big\{\sqrt{\pi_{g}(x_i)\pi_{f_0}(x_i)}+\sqrt{(1-\pi_{g}(x_i))(1-\pi_{f_0}(x_i))}\Big\}.
\end{eqnarray*}
In other words we have
\begin{multline*}
\log\mathbb{E}\Big[\exp\Big(\frac{n}{2}\Big(\gamma_n(f_{0})-\gamma_n(g)\Big)\Big)=\\
\sum_{i=1}^{n}\log\Big\{ 1-\frac{1}{2}\Big[\Big(\sqrt{\pi_{f_0}(x_i)}-\sqrt{\pi_{g}(x_i)}\Big)^{2} +\Big(\sqrt{1-\pi_{f_0}(x_i)}-\sqrt{1-\pi_{g}(x_i)}\Big)^{2}\Big]\Big\}.
\end{multline*}
This implies that
\begin{eqnarray*}
\log\mathbb{E}\Big[\exp\Big(\frac{n}{2}\Big(\gamma_n(f_{0})-\gamma_n(g)\Big)\Big)\Big]
&\leq& \sum_{i=1}^{n}-\frac{1}{2}\Big[\Big(\sqrt{\pi_{f_0}(x_i)}-\sqrt{\pi_{g}(x_i)}\Big)^{2} +\Big(\sqrt{1-\pi_{f_0}(x_i)}-\sqrt{1-\pi_{g}(x_i)}\Big)^{2}\Big]\\
&=& -nh^{2}(\mathbb{P}_{f_0},\mathbb{P}_{g}).
\end{eqnarray*}
Consequently
$$\mathbb{P}(\overline{\gamma}_n(f_{0})-\overline{\gamma}_n(g)\geq b)\leqslant\exp\Big[\frac{-nb}{2}+ \frac{n}{2}\mathcal{K}(\mathbb{P}_{f_0}^{(n)},\mathbb{P}_{g}^{(n)})-nh^{2}(\mathbb{P}_{f_0}^{(n)},\mathbb{P}_{g}^{(n)})\Big],$$
and, if we choose for positive $x$,
$$b= \frac{2x}{n}+\mathcal{K}(\mathbb{P}_{f_0}^{(n)},\mathbb{P}_{g}^{(n)})-2h^{2}(\mathbb{P}_{f_0}^{(n)},\mathbb{P}_{g}^{(n)})>0,$$
we  have,
$$\mathbb{P}\Big(\overline{\gamma}_n(f_{0})-\overline{\gamma}_n(g)\geq \frac{2x}{n}+\mathcal{K}(\mathbb{P}_{f_0}^{(n)},\mathbb{P}_{g}^{(n)})-2h^{2}(\mathbb{P}_{f_0}^{(n)},\mathbb{P}_{g}^{(n)})\Big)\leqslant\exp(-x).$$
We conclude that
$ \mathbb{P}(\Omega_2^{c}(\xi))\leqslant\sum_{m^{\prime}}\exp(-L_m^{\prime}\vert m^{\prime}\vert-\xi)\leq\Sigma e^{-\xi},$
which ends the  proof of (R-\ref{biais2}).\qed
\section{Appendix}\label{S6}
\setcounter{equation}{0}
\setcounter{lem}{0}
\setcounter{theo}{0}
\label{appendix}

\subsection{Proof  of Lemma \ref{Projhisto}.}
By definition
$$f_m=\arg\min_{f \in S_m} \left[ \sum_{i=1}^n \log(1+\exp(f(x_i)))-\pi_{f_0}(x_i)f(x_i)   \right].$$
For all $f\in S_m$, for all $J\in m$ and for all $x\in J$, we have $f(x)=f^{(J)}$. Hence $f_m(x)=\overline{f}_m^{(J)}$ for all $x$ in $J$, and
for all $J$ in $m$, we aim at finding $\overline{f}_m^{(J)}$ such that
$$\overline{f}_m^{(J)}=\arg\min_{f^{(J)} } \left[ |J|\log(1+\exp(f^{(J)}))-\sum_{i\in J}\pi_{f_0}(x_i)f^{(J)}   \right]$$
where $|J|={\mbox card} \{i\in\{1,...,n\} ; x_i\in J\}$. 
Easy calculations show that he coefficient $\overline{f}^{(J)}_m$ satisfies
$$\vert J\vert \frac{\exp(\overline{f}_m^{(J)})}{1+\exp (\overline{f}_m^{(J)})}-\sum_{i\in J}\pi_{f_0}(x_i)=0,$$
that is
\begin{eqnarray}
\label{projhisto}\overline{f}_m^{(J)}=\log \left(  
\frac{ \sum_{i\in J}\pi_{f_0}(x_i)  }{\vert  J \vert  (1-\sum_{i\in J}\pi_{f_0}(x_i) / \vert  J \vert ) }\right).
\end{eqnarray}
Consequently, 
$\pi_{f_m}$ defined as in \eref{pif}
satisfies
that $\pi_{f_m}(x)=\pi_{f_m}^{(J)}$ for all $x\in J$, where
$$\pi_{f_m}^{(J)}=\frac{1}{\vert J\vert}\sum_{i\in J}\pi_{f_0}(x_i),$$
and hence $\pi_{f_m}=\arg\min_{t\in S_m}\parallel t-\pi_{f_0}\parallel_n$
is the usual projection of $\pi_{f_0}$ on to $S_m=<\Phi_j, j\in m>.$
In the same way, $\hat{f}_m$ defined by \eref{fchapD}
satisfies $\hat f_m(t)=\hat{f}_m^{(J)} $ for all $t\in J$, where
$$\hat{f}_m^{(J)}=\log \left(  
\frac{ \sum_{i\in J}Y_i  }{\vert  J \vert  (1-\sum_{i\in J}Y_i/ \vert  J \vert ) }\right).$$ In other words,
$\pi_{\hat f_m}$, defined as $\pi_f$ with $f$ replaced by $\pi_{\hat f_m}$,
satisfies $\pi_{\hat f_m}(x)=\pi_{\hat f_m}^{(J)}$, \mbox{ for all} $ x\in J$, with $$\pi_{\hat f_m}^{(J)}=\frac{1}{\vert J\vert}\sum_{i\in J}Y_i.$$

\subsection{Proof  of Lemma \ref{lm}.}

In the following, for the sake of notation simplicity,  we will use   $\gamma(\beta)$ for $\gamma(f_{\beta})$.
A second-order Taylor expansion of the function $\gamma()$ around $\beta^{*}$ gives for any $\beta\in \Lambda_m$
\begin{multline*}
\gamma(\beta)=\gamma(\beta^{*})+\nabla_{\beta}\gamma(\beta^{*})(\beta-\beta^{*})
\\+ \int_{0}^{1}(1-t)\sum_{i_1+\dots+i_D=2}\frac{2!}{i_1!\dots i_D!}(\beta_1-\beta_{1}^{*})^{i_1}\dots 
(\beta_D-\beta_{D}^{*})^{i_D}\frac{\partial\gamma^{2}}{\partial\beta_1\dots\partial\beta_D}(\beta^{*}+t(\beta-\beta{*}))dt .
\end{multline*}
 Easy calculation shows that 
\begin{eqnarray*}
&& \sum_{i_1+\dots+i_D=2}\frac{2!}{i_1!\dots i_D!}(\beta_1-\beta_{1}^{*})^{i_1}\dots 
(\beta_D-\beta_{D}^{*})^{i_D}\frac{\partial\gamma^{2}}{\partial\beta_1\dots\partial\beta_D}(\beta^{*}+t(\beta-\beta{*}))\\
&=&\sum_{j=1}^{D}\frac{1}{n}\sum_{i=1}^{n}\psi_{j}^{2}(x_{i})(\beta_{j}-\beta^{*}_{j})^{2}\pi\left(f_{\beta^{*}+t(\beta-\beta{*})}(x_i)\right)\left[1-\pi\left(f_{\beta^{*}+t(\beta-\beta{*})}(x_i)\right)\right]\\
&+&2\sum_{l\neq k}\frac{1}{n}\sum_{i=1}^{n}\psi_{l}(x_{i})\psi_{k}(x_i)(\beta_l-\beta_l^{*})(\beta_k-\beta_k^{*})\pi\left(f_{\beta^{*}+t(\beta-\beta{*})}(x_i)\right)\left[1-\pi\left(f_{\beta^{*}+t(\beta-\beta{*})}(x_i)\right)\right]\\
&=&\frac{1}{n}\sum_{i=1}^{n}\pi\left(f_{\beta^{*}+t(\beta-\beta{*})}(x_i)\right)\left[1-\pi\left(f_{\beta^{*}+t(\beta-\beta{*})}(x_i)\right)\right](f_\beta(x_i)-f_{\beta^{*}}(x_i))^{2}.
\end{eqnarray*}
 This implies that
$$\gamma(\beta)\geq\gamma(\beta^{*})+\nabla_\beta\gamma(\beta^{*})(\beta-\beta^{*})+\frac{\mathcal{U}_{0}^{2}}{2}\lVert f_{\beta}-f_{\beta^{*}}\rVert^2_{n}.$$
Since $\beta^{*}$ is the minimizer of $\gamma(.)$ over the set $\Lambda_m$, we have $\nabla_\beta\gamma(\beta^{*})(\beta-\beta^{*})\geq 0$ for all
$\beta \in \Lambda_m$. Thus the result follows.

\subsection{Proof of Lemma~\ref{control}}
Let $S_{D}$ and $S_{D^{\prime}}$ two vector spaces of dimension $D$ and $D^{\prime}$ respectively. Set $S=S_{D}\cap \mathbb{L}_\infty(C_0)+S_{D^{\prime}}\cap \mathbb{L}_\infty(C_0)$  and $\vec{\varepsilon}^{\prime}$ be an independent copie of $\vec{\varepsilon}.$ Set 
\begin{equation}\label{Z}
 Z=\sup_{u\in S}\frac{\langle\vec{\varepsilon} ,u\rangle_n}{\parallel u \parallel_n} ,
\mbox{ and for all }i=1,\dots,n, \quad Z^{(i)}=\sup_{u\in S}\frac{1}{\parallel u \parallel_n}\left(\frac{1}{n}\sum_{k\neq i}\varepsilon_{k}u(x_{k})+ \varepsilon_{i}^{\prime}u(x_{i}) \right).
\end{equation}
By Cauchy-Schwarz Inequality the supremum in (\ref{Z}) is achieved at $\Pi_S(\vec{\varepsilon}).$ Consequently,
$$Z-Z^{(i)}\leqslant \frac{(\varepsilon_{i}-\varepsilon_{i}^{\prime})(\Pi_S(\vec{\varepsilon})(x_i)}{n\parallel\Pi_S(\vec{\varepsilon})\parallel_n},
 \qquad \mbox{ and }\qquad 
 \mathbb{E}_{f_0}[(Z-Z^{(i)})^{2}|\vec{\varepsilon}]
\leq
\mathbb{E}_{f_0}\left[\frac{(\varepsilon_{i}-\varepsilon_{i}^{\prime})^{2}[\Pi_S(\vec{\varepsilon})(x_i)]^2}{n^2\parallel\Pi_S(\vec{\varepsilon})\parallel_n^2}|\vec{\varepsilon}\right]
$$
with 
\begin{eqnarray*}
\mathbb{E}_{f_0}\left[\frac{(\varepsilon_{i}-\varepsilon_{i}^{\prime})^{2}[\Pi_S(\vec{\varepsilon})(x_i)]^2}{n^2\parallel\Pi_S(\vec{\varepsilon})\parallel_n^2}|\vec{\varepsilon}\right]
&=&\frac{[\Pi_S(\vec{\varepsilon})(x_i)]^2}{n^2\parallel\Pi_S(\vec{\varepsilon})\parallel_n^2}\mathbb{E}_{f_0}\left[(\varepsilon_{i}-\varepsilon_{i}^{\prime})^{2}|\vec{\varepsilon} \right]\\&=&\frac{[\Pi_S(\vec{\varepsilon})(x_i)]^2}{n^2\parallel\Pi_S(\vec{\varepsilon})\parallel_n^2}\left(    \varepsilon_i^2+\mathbb{E}_{f_0}(\varepsilon_i^2) \right)\leq \frac{5[\Pi_S(\vec{\varepsilon})(x_i)]^2}{4n^2\parallel\Pi_S(\vec{\varepsilon})\parallel_n^2} .
\end{eqnarray*}
This implies that
$$
\sum_{i=1}^{n}\mathbb{E}_{f_0}[(Z-Z^{(i)})^{2}\ind_{Z>Z^{(i)}}|\vec{\varepsilon}]
\leq \frac{5}{4n}.
$$
We now apply Lemma \ref{Bouch}  from  Boucheron \textit{et al.} \citeyear{boucheron}), that is recalled here.
\begin{lem}\label{Bouch} 
Let $X_{1} ,\dots, X_{n}$ independent random variables taking values in a measurable space $\mathcal{X}$.
 Denote by $X_{1}^{n}$ the vector of these $n$ random variables. Set $Z=f(X_{1},\dots,X_{n})$ ~ and ~ $Z^{(i)}=f(X_{1},\dots,X_{i-1},X_{i}^{\prime},X_{i+1},\dots,X_{n}),$
where $X_{1}^{\prime} ,\dots, X_{n}^{\prime}$ denote independent copies of $X_{1} ,\dots, X_{n}$ and f : $\mathcal{X}^{n} \rightarrow \mathbb{R}$ some measurable function.
Assume that there exists a positive constant $c$ such that,
$\mathbb{E}_{f_0}\left[\sum_{i=1}^{n}(Z-Z^{(i)})^{2}\mathds{1}_{Z> Z^{(i)}}|X_{1}^{n}\right]\leqslant c$. Then for all $t > 0$,
$$\mathbb{P}_{f_0} (Z>\mathbb{E}_{f_0}(Z)+t)\leqslant e^{-t^{2}/4c}.$$
\end{lem}
Applying Lemma~\ref{Bouch} to $Z$ defined in \eref{Z}, we obtain that for all $x>0$,
$$\mathbb{P}\left(\sup_{u\in S}\frac{\langle\vec{\varepsilon} ,u\rangle_n}{\parallel u \parallel_n}>
\mathbb{E}_{f_0}\left[\sup_{u\in S}\frac{\langle\vec{\varepsilon} ,u\rangle_n}{\parallel u \parallel_n}\right]+\sqrt{\frac{5x}{n}}\right)\leqslant \exp{(-x)}.$$
Let $\{ \psi_1,\dots,\psi_{D+D^{\prime}}\}$ be an orthonormal basis of $S_{D}+S_{D^{\prime}}$. Using Jensen's Inequality, we write
\begin{align*}
 \mathbb{E}_{f_0}\left[\sup_{u\in S}\frac{\langle\vec{\varepsilon} ,u\rangle_n}{\parallel u \parallel_n}\right]=
\mathbb{E}_{f_0}(\parallel\Pi_S(\vec{\varepsilon})\parallel_n)&=\mathbb{E}_{f_0}\left[\left(\sum_{k=1}^{D+D^{\prime}}(\langle\vec{\varepsilon} ,\psi_k\rangle_n)^2\right)^{1/2}\right]\\
&\leq\left(\sum_{k=1}^{D+D^{\prime}}\mathbb{E}_{f_0}(\langle\vec{\varepsilon} ,\psi_k\rangle_n)^2\right)^{1/2}\\
&\leqslant \sqrt{\frac{ D+D^{\prime}}{4n}}.
\end{align*}
This concludes the proof of Lemma~\ref{control}.


\bibliography{biblioM}
 \bibliographystyle{chicago}

\end{document}